\documentclass[12pt,a4paper,draft]{article}%
\usepackage{graphicx}
\usepackage{amsmath}
\usepackage{amsfonts}
\usepackage{amssymb}%
\providecommand{\U}[1]{\protect\rule{.1in}{.1in}}
\providecommand{\U}[1]{\protect\rule{.1in}{.1in}}
\newtheorem{theorem}{Theorem}

\newtheorem{corollary}[theorem]{Corollary}

\newtheorem{definition}[theorem]{Definition}

\newtheorem{lemma}[theorem]{Lemma}

\newtheorem{proposition}[theorem]{Proposition}
\newtheorem{remark}[theorem]{Remark}

\newenvironment{proof}[1][Proof]{\noindent\textbf{#1.} }{\ \rule{0.5em}{0.5em}}
\begin{document}

\title{On a hierarchy of infinite-dimensional spaces and related Kolmogorov-Gelfand widths}
\author{O. Kounchev\\Institute of Mathematics and Informatics, \\Bulgarian Academy of Sciences\\and\\IZKS, University of Bonn}
\maketitle

\begin{abstract}
Recently the theory of widths of Kolmogorov-Gelfand has received a great deal
of interest due to its close relationship with the newly born area of
Compressed Sensing. It has been realized that widths reflect properly the
sparsity of the data in Signal Processing. However fundamental problems of the
theory of widths in multidimensional Theory of Functions remain untouched, as
well as analogous problems in the theory of multidimensional Signal Analysis.
In the present paper we provide a multidimensional generalization of the
original result of Kolmogorov by introducing a new hierarchy of
infinite-dimensional spaces based on solutions of higher order elliptic equation.

\end{abstract}

\section{Introduction}

Recent interest to the theory of widths (especially to Gelfand widths) has
been motivated by applications in Compressed Sensing (CS). In a certain sense
the central idea of CS is rooted in the theory of widths, cf. e.g.
\cite{donoho}, \cite{candes}, \cite{devore}, \cite{pinkus2011}. However,
apparently this strategy works smoothly only in the case of representation of
one-dimensional signals, while an adequate approach to multivariate signals is
missing -- one reason may be found by analogy in the fact that the theory of
Kolmogorov-Gelfand widths fits properly only for one-dimensional function
spaces (as pointed out below, e.g. in formula (\ref{dNisInfinity})). Recently,
a new multivariate Wavelet Analysis was developed based on solutions of
elliptic partial differential equations (\cite{okbook}), in particular
"polyharmonic subdivision wavelets" were introduced (cf.
\cite{dynkounchevlevinrender}, \cite{kounchevKalaglarsky}); in order to apply
CS ideas to these wavelets it would require essential generalization of the
theory of widths for infinite-dimensional spaces. We start with this
motivation to study a generalization of the Kolmogorov-Gelfand theory of
widths by introducing a \emph{new hierarchy of infinite-dimensional spaces}
based on solutions of higher order elliptic equations. However, there is a
different perspective on the present research: its main purpose is to
introduce this new hierarchy and to apply it to the theory of widths as a
testing field. One may expect also that this development would throw a new
light on the nature of \emph{sparsity}\textbf{ }in multidimensional Signal Analysis.

In his seminal paper \cite{kolmogorov1936} Kolmogorov has introduced the
theory of widths and has applied it ingeniously to the following set of
functions defined in the compact interval:
\begin{equation}
K_{p}:=\left\{  f\in AC^{p-1}\left(  \left[  a,b\right]  \right)  :%
{\displaystyle\int_{0}^{1}}
\left\vert f^{\left(  p\right)  }\left(  t\right)  \right\vert ^{2}%
dt\leq1\right\}  . \label{Kp}%
\end{equation}
In the present paper we study a natural multivariate generalization of the set
$K_{p}$ which in a domain $B\subset\mathbb{R}^{n}$ is given by
\begin{equation}
K_{p}^{\ast}:=\left\{  u\in H^{2p}\left(  B\right)  :%
{\displaystyle\int_{B}}
\left\vert \Delta^{p}u\left(  x\right)  \right\vert ^{2}dx\leq1\right\}  ,
\label{Kpstar}%
\end{equation}
where $\Delta^{p}$ is the $p-$th iterate of the Laplace operator $\Delta=%
{\displaystyle\sum_{j=1}^{n}}
\partial^{2}/\partial x_{j}^{2};$ we consider more general sets $K_{p}^{\ast}$
given in (\ref{KpstarGeneral}) below. We generalize the notion of width in the
framework of the Polyharmonic Paradigm, and obtain analogs to the
one-dimensional results of Kolmogorov.

The \emph{Polyharmonic Paradigm} has been announced in \cite{okbook} as a new
approach to Multidimensional Mathematical Analysis, which is based on
solutions of higher order elliptic partial differential equations as opposed
to the usual concept which is based on algebraic and trigonometric polynomials
of several variables. It has proved to be very successful in the Moment
Problem \cite{kounchevrenderMoment}, Approximation and Interpolation
\cite{kounchev1991Hanstholm}, \cite{kounchev1992TAMS}, and Spline Theory
\cite{kounchev1998TAMS}, \cite{okbook}.

The main objective of the present research is a new development of the
Polyharmonic Paradigm. It provides a new hierarchy of infinite-dimensional
spaces of functions which are used for a generalization of the Kolmogorov's
theory of widths. This new hierarchy generalizes the usual hierarchy of
finite-dimensional subspaces $X_{N}$ of the space $C^{\infty}\left(  I\right)
$ for an interval $I\subset\mathbb{R}.$ The crux of this notion of hierarchy
is the following: Let the domain $D\subset\mathbb{R}^{n},$ be compact with
sufficiently smooth boundary $\partial D.$ Then the $N-$dimensional subspaces
in $C^{\infty}\left(  I\right)  $ will be generalized by spaces of
\emph{solutions of elliptic equations} (and by more general spaces introduced
in Definition \ref{Dhdimension} below):
\begin{equation}
X_{N}=\left\{  u:P_{2N}u\left(  x\right)  =0,\quad\text{for }x\in D\right\}
\subset L_{2}\left(  D\right)  ; \label{XN}%
\end{equation}
here $P_{2N}$ is an elliptic operator of order $2N$ in the domain $D.$
Respectively, the simplest version of our generalization of Kolmogorov's
theorem about widths finds the extremizer of the following problem
\[
\inf_{X_{N}}\operatorname*{dist}\left(  X_{N},K_{p}^{\ast}\right)  ,
\]
where $K_{p}^{\ast}$ is the set defined in (\ref{Kpstar}) and $X_{N}$ is
defined in (\ref{XN}) for arbitrary elliptic operator $P_{2N}$ of order $2N;$
for the complete formulation see Theorem \ref{TKolmogorovMultivariate} below.

What is the reason to take namely solutions of elliptic equations in the
multidimensional case is explained in the following section.

\subsection{The hierarchy of infinite-dimensional spaces - a justification via
Chebyshev systems\label{SIntroHierarchy}}

Let us give a heuristic outline of the \emph{motivation} and \emph{the main
idea} of this new hierarchy of spaces, by explaining how it appears as a
natural generalization of the finite-dimensional subspaces of $C^{N-1}\left(
\left[  a,b\right]  \right)  $ in a compact interval $\left[  a,b\right]  $ in
$\mathbb{R}.$

First of all, let us understand the structure of the finite-dimensional
subspaces of $C^{N-1}\left(  \left[  a,b\right]  \right)  $: It is important
to note that for a "general position" (or "generic") $N-$dimensional subspace
$X_{N}\subset C^{N-1}\left(  I\right)  $ in the interval $I=\left(
a,b\right)  ,$ there exists a \emph{finite or infinite number} of subintervals
$I_{k,j}=\left(  a_{k,j},b_{k,j}\right)  $ with $\bigcup\overline{I_{k,j}%
}=\left[  a,b\right]  ,$ and a basis $\left\{  u_{k}\right\}  _{k=1}^{N},$
\begin{equation}
X_{N}=\operatorname*{span}\left\{  u_{k}\right\}  _{k=1}^{N},\label{XN1dim}%
\end{equation}
(here $\operatorname*{span}$ denotes the linear closure) where the Wronski
determinants satisfy
\begin{align}
\varepsilon_{k,j}W\left(  u_{1}\left(  t\right)  ,u_{2}\left(  t\right)
,...,u_{k}\left(  t\right)  \right)   &  >0\qquad\text{for }t\in
I_{k,j}\label{Wronski}\\
\text{with }\varepsilon_{k,j} &  =1\text{ or }-1.\label{Wronski+}%
\end{align}
A simplest example would be the space $X_{2}:=\operatorname*{span}\left\{
1,t^{2}\right\}  $ considered on the interval $\left[  -1,1\right]  ,$ where
the Wronskian $W\left(  1,t^{2}\right)  $ changes sign at $0.$\footnote{Note
that there are cases where the system of functions $\left\{  u_{k}\right\}
_{k=1}^{N}$ has dimension $N$ but its Wronskian is $0$ on a whole interval,
e.g. the system of two functions $\left\{  t^{2},\chi\left(  t\right)
t^{2}-\chi\left(  -t\right)  t^{2}\right\}  $ in the interval $\left[
-1,1\right]  ,$ where $\chi$ is the Heaviside function. However this is an
exception, hence not "generic". One may try to make the last precise:\ By
introducing a proper topology/metric in the set $\mathcal{S}_{N}$ of all
$N-$dimensional spaces in $C^{N}\left(  I\right)  ,$ e.g. by taking the
distance between the unit spheres in two spaces, we may prove that in the set
$\mathcal{S}_{N}$ those spaces having Wronskian equal to $0$ are a "small set"
in the sense of second category of Baire.}

Since A. Markov it is known that the positivity condition (\ref{Wronski}%
)-(\ref{Wronski+}) is characteristic for Extended Complete Chebyshev systems
(called  $ECT-$systems in \cite{karlinstudden}, chapter $11;$  cf. also
\cite{kreinnudelman}, chapter $2,$ section $5,$ Theorem $5.1$). For that
reason, we may formulate our \textbf{important observation} by saying  that a
"general position" $N-$dimensional space $X_{N}\subset C^{N}\left(  I\right)
$ is a \emph{piecewise Extended Complete Chebyshev system} of order $N.$ 

We will remind some basic properties related to Extended Complete Chebyshev
Systems. The following fundamental result describes their structure  (cf.
\cite{karlinstudden}, chapter $11,$ Theorem $1.1$).

\begin{proposition}
\label{PstructureMarkovSystem} Let us assume that the space $X_{N}$ in
(\ref{XN1dim}) restricted to $\ $some subinterval $J\subset I$ has
sign-definite Wronskians as in (\ref{Wronski})-(\ref{Wronski+}). Then the
restriction of the space $X_{N}$ to the interval  $J$ is a space of solutions
for an ODE
\begin{equation}
L_{N}^{J}u\left(  t\right)  =0,\qquad\text{for }t\in J;\label{LN+1}%
\end{equation}
here $L_{N}^{J}$ is an ordinary differential operator of order $N$ in $J,$
given by
\begin{equation}
L_{N}^{J}\left(  t;\frac{d}{dt}\right)  =%
{\displaystyle\prod_{k=1}^{N}}
\frac{d}{dt}\frac{1}{\rho_{k}\left(  t\right)  }\label{LNoperatorForm}%
\end{equation}
where the functions $\rho_{k}$ satisfy $\rho_{k}\left(  t\right)  >0$ in $J.$
\end{proposition}

Let us remark that the weight functions $\rho_{k}$ may be chosen in different
ways, cf. \cite{mazure}. If we put for the Wronskians
\[
W_{k}:=W\left(  u_{1}\left(  t\right)  ,u_{2}\left(  t\right)  ,...,u_{k}%
\left(  t\right)  \right)
\]
then the functions $\rho_{k}\left(  t\right)  $ may be written  as%
\begin{align*}
\rho_{1} &  =W_{1}=u_{1},\quad\rho_{2}=W_{2}/W_{1}^{2}\\
\rho_{k} &  =W_{k}W_{k-2}/W_{k-1}^{2}\qquad\text{for }k\geq3,
\end{align*}
(cf. \cite{polya},  \cite{kreinnudelman}, section $5,$ chapter $2,$ Theorem
$5.1,$ or \cite{karlinstudden}, chapter $11$).

Obviously, the operator $L_{N}$ has a non-negative leading coefficient and is
in this sense one-dimensional "elliptic".

\begin{remark}
The detailed proof of Proposition \ref{PstructureMarkovSystem} is a part of
the general theory of Chebyshev systems developed by A. Markov,  S. Bernstein,
M. Krein and others,  in which the Extended Complete Chebyshev systems are a
special case which are of interest for us. Their theory is presented in detail
in the above mentioned monographs \cite{kreinnudelman} and
\cite{karlinstudden}, whereby in the first reference the case of
non-differentiable systems is emphasized.

Let us mention that the space $X_{N}$ generated by a Chebyshev system is often
called \textbf{Haar space, }cf. \cite{mccullough}. Thus one may also say that
a generic $N-$dimensional subspace of $C^{N-1}\left(  I\right)  $ is piecewise
Haar space.
\end{remark}

First, we will seek for a \emph{generalizable framework} for the Extended
Complete Chebyshev systems which we have obtained on every subinterval
$I_{j}.$ The interpolation framework seems to be suitable: Let us note that
condition (\ref{Wronski}) has equivalent formulation as \emph{Hermite
interpolation},  in particular, for arbitrary $t_{0}\in\left(  a,b\right)  $
and constants $\left\{  c_{k}\right\}  _{k=0}^{N-1},$ it is possible to solve
the interpolation problem
\begin{equation}
u^{\left(  k\right)  }\left(  t_{0}\right)  =c_{k}\qquad\text{for
}k=0,1,...,N-1\label{Dirichlet}%
\end{equation}
where $u\in X_{N}.$

At this point it is important to emphasize that we will select judiciously,
and generalize in the multivariate case, \emph{only} those $2M-$dimensional
subspaces $X_{2M}=\operatorname*{span}\left\{  u_{k}\right\}  _{k=1}%
^{2M}\subset C^{2M-1}\left(  I\right)  $ which satisfy a rather specific
interpolation property:

\begin{definition}
\label{DDirichletBVP}We say that the space $X_{2M}\subset C^{2M-1}\left(
I\right)  $ has the \textbf{Dirichlet BVP property}, if for every subinterval
$I_{1}=\left[  a_{1},b_{1}\right]  \subset I,$ and for arbitrary constants
$\left\{  c_{k},d_{k}\right\}  _{k=0}^{M-1},$ the (Dirichlet) boundary value
problem
\begin{align}
u^{\left(  k\right)  }\left(  a_{1}\right)   &  =c_{k}\qquad\text{for
}k=0,1,...,M-1\label{Dirichlet1}\\
u^{\left(  k\right)  }\left(  b_{1}\right)   &  =d_{k}\qquad\text{for
}k=0,1,...,M-1 \label{Dirichlet2}%
\end{align}
has a solution $u\in X_{2M}.$
\end{definition}

\begin{remark}
Let us assume  that in the space $X_{2M}$ there exists an Extended Complete
Chebyshev system in $I$ (i.e. a system satisfying positivity of the Wronskians
(\ref{Wronski})-(\ref{Wronski+}) in $I$). Then  $X_{2M}$ satisfies Definition
\ref{DDirichletBVP} which follows from the very  definition of Extended
Complete Chebyshev systems, cf. \cite{karlinstudden}, chapter $11$. Thus the
Extended Complete Chebyshev systems provide the main bulk of examples for
Definition \ref{DDirichletBVP}.
\end{remark}

One may consider the solvability of problem (\ref{Dirichlet1}%
)-(\ref{Dirichlet2}) as a "parametrization" of the space $X_{2M}$ by the
Dirichlet boundary values $\left\{  c_{k},d_{k}\right\}  _{k=0}^{M-1},$ and
this important property will be generalized to the multivariate case.

We are interested in the BVP interpretation which follows from Proposition
\ref{PstructureMarkovSystem}: Since the space $X_{2M}$ may be represented as
\[
X_{2M}=\left\{  u:L_{2M}u=0\quad\text{for }t\in I\right\}  ,
\]
for an elliptic operator $L_{2M},$ then the solvability of problem
(\ref{Dirichlet1})-(\ref{Dirichlet2}) in the space $X_{2M}$ may be considered
as a special case of the multidimensional theory for \textbf{Elliptic Boundary
Value Problems} (\textbf{BVP}), and it is a classical BVP in the
one-dimensional ODEs as well, cf. \cite{naimark}.

In view of the last observation, we seek for a multidimensional generalization
of problem (\ref{Dirichlet1})-(\ref{Dirichlet2}). Let $D$ be a bounded domain
in $\mathbb{R}^{n}$ and consider the subspaces of $L_{2}\left(  D\right)  .$
The space of solutions of an elliptic equation generalizing equation
(\ref{LN+1}) may be considered as a natural generalization of the space
$X_{2M}.$ Indeed, if
\begin{equation}
X_{2M}=\left\{  u:P_{2M}u\left(  x\right)  =0\quad\text{for }x\in D\right\}
,\label{S2M}%
\end{equation}
where $P_{2M}\left(  x;D_{x}\right)  $ is an elliptic differential operator in
the domain $D,$ then the natural generalization to problem (\ref{Dirichlet1}%
)-(\ref{Dirichlet2}) is an Elliptic BVP, as for example the Dirichlet problem
which may be considered for subdomains $D_{1}$ in $D,$ namely
\begin{align}
P_{2M}u\left(  x\right)   &  =0\quad\text{for }x\in D_{1}%
\label{DirichletMultivariate}\\
\left(  \frac{\partial}{\partial n}\right)  ^{k}u\left(  y\right)   &
=c_{k}\left(  y\right)  \qquad\text{for }y\in\partial D_{1},\quad\text{for
}k=0,1,...,M-1.\label{DirichletMultivariate2}%
\end{align}
Let us remind that the Dirichlet problem is well-known to be solvable for data
$\left\{  c_{k}\left(  y\right)  \right\}  _{k=0}^{M-1}$ from a proper Sobolev
or H\"{o}lder space on the boundary $\partial D_{1}.$ An important point is
that for  a large class of operators $P_{2M}$ every solution of
(\ref{DirichletMultivariate})-(\ref{DirichletMultivariate2}) may be
approximated by solutions in the whole domain $D,$ i.e. by elements of
$X_{2M}.$ This may be considered as a substitute of the interpolation property
(\ref{Dirichlet1})-(\ref{Dirichlet2}) in the one-dimensional case. Very
important hint for identifying the operators $P_{2M}$ which represent
Multidimensional Chebyshev systems is provided by formula
(\ref{LNoperatorForm}). This is the main reason for the judicious choice of
the special class of operators $P_{2M}$ in Definition \ref{Dhdimension} below,
as they mimic the operators in (\ref{LNoperatorForm}) and serve our
purposes.\footnote{This generalization has been discussed in detail in
\cite{kounchev2008-ChebyshevSystems}.}

Making analogy with the one-dimensional case (\ref{Dirichlet1}%
)-(\ref{Dirichlet2}), we may say that here the space $X_{2M}$ defined in
(\ref{S2M}) is "parametrized" by the Dirichlet boundary conditions
(\ref{DirichletMultivariate}), however the "parameter" $\left\{  c_{k}\left(
y\right)  \right\}  _{k=0}^{M-1}$ runs a function space. Hence, the spaces
$X_{2M}$ may be considered as a natural generalization of the one-dimensional
Extended Complete Chebyshev systems and we call them \textbf{Multidimensional
Chebyshev systems}.

After we have the Multidimensional Chebyshev systems in our disposal, the next
step will be to introduce the multivariate generalization of the
$N-$dimensional subspaces of $C^{\infty}\left(  I\right)  .$ We will define
them in Definition \ref{Dhdimension} below as subspaces $X_{N}$ of functions
in $C^{\infty}\left(  D\right)  $ which are piecewise solutions of (regular)
elliptic differential operators of order $2N.$ We will say that $X_{N}$ has
"Harmonic Dimension $N$" and we will write
\[
\operatorname{hdim}\left(  X_{N}\right)  =N,
\]
see Definition \ref{Dhdimension} below. Kolmogorov's notion of $N-$width (and
in a similar way Gelfand's width) is naturally generalized for symmetric sets
by the notion of "Harmonic $N-$width" defined by putting
\[
\operatorname*{hd}\nolimits_{N}\left(  S\right)  :=\inf_{\operatorname{hdim}%
\left(  X_{N}\right)  =N}\operatorname*{dist}\left(  X_{N},S\right)  ,
\]
see Definition \ref{Dwidth} below. The \textbf{main result} of the present
paper is the computation of
\[
\operatorname*{hd}\nolimits_{N}\left(  K_{p}^{\ast}\right)  \qquad\text{for
}N\leq p,
\]
where $K_{p}^{\ast}$ is defined in (\ref{Kpstar}) and more generally in
(\ref{KpstarGeneral}).

\subsection{Plan of the paper}

To facilitate the reader, in section \ref{SKolmogorovsResults} we provide a
short summary of the original Kolmogorov's results. For the same reason, in
section \ref{SreminderBVP} we provide a short reminder on Elliptic BVP. In
section \ref{SprincipalAxes} we prove the representation of the "cylindrical
ellipsoid" set $K_{p}^{\ast}$ in principal axes which generalizes the
one-dimensional representation of Kolmogorov, cf. Theorem
\ref{TprincipalAxesMultivariate} below. In section \ref{Shierarchy} we
introduce the notion of \emph{Harmonic Dimension}, and the \emph{First Kind}
spaces of Harmonic Dimension $N.$ Based on it we define \emph{Harmonic Widths}
which generalize Kolmogorov's widths. In section \ref{Swidths}, in Theorem
\ref{TKolmogorovMultivariate} we prove a genuine analog to Kolmogorov's
theorem about widths. It says that among all spaces $X_{N}$ having Harmonic
Dimension $N,$ some special space $\widetilde{X}_{N}$ provides the best
approximation to the set $K_{p}^{\ast}$ in problem
\[
\inf_{X_{N}}\operatorname*{dist}\left(  X_{N},K_{p}^{\ast}\right)  ,
\]
and this space $\widetilde{X}_{N}$ is identified by the principal axes
representation provided by Theorem \ref{TprincipalAxesMultivariate}. In
section \ref{SsecondKind} we introduce \emph{Second Kind} spaces of Harmonic
Dimension $N$ and formulate a further generalization of Theorem
\ref{TKolmogorovMultivariate}. Apparently, the First and Second Kind spaces
having Harmonic Dimension $N$ provide the maximal generalization in the
present framework.

A special case of the present results is available in \cite{kounchevSozopol},
and might be instructive for the reader to start with.

A final remark to our generalization is in order. In our consideration we will
not strive to achieve a maximal generality. As it is clear, especially in the
applications to the theory of widths even in the one-dimensional case we may
consider not all $N-$dimensional subspaces but "almost all" $N-$dimensional
subspaces of $C^{\infty}\left(  D\right)  $ in some sense, or a class of
$N-$dimensional subspaces which are dense (in a proper topology) in the set of
all other $N-$dimensional subspaces. This "genericity" point of view is
essential in our multivariate generalization since it will allow us to avoid
burdensome proofs necessary in the case of the bigger generality of the
construction. For the same reason we will not consider elliptic
pseudo-differential operators although almost all results have a
generalization for such setting.

\textbf{Acknowledgements:} The author acknowledges the support of the
Alexander von Humboldt Foundation, and of Project DO-02-275 with Bulgarian
NSF. The author thanks the following Professors: Matthias Lesch for the
interesting discussion about hierarchies of infinite-dimensional linear
spaces, Hermann Render about advice on multivariate polynomial division, and
Peter Popivanov, Nikolay Kutev and Georgi Boyadzhiev about advice on Elliptic BVP.

\section{Kolmogorov's results - a reminder \label{SKolmogorovsResults}}

In order to make our multivariate generalization transparent we will recall
the original results of Kolmogorov provided in his seminal paper
\cite{kolmogorov1936}. Kolmogorov has considered the set $K_{p}$ defined in
(\ref{Kp}). He proved that this is an \textbf{ellipsoid} by constructing
explicitly its principal axes. Namely, he considered the eigenvalue problem
\begin{align}
\left(  -1\right)  ^{p}u^{\left(  2p\right)  }\left(  t\right)   &  =\lambda
u\left(  t\right)  \qquad\qquad\qquad\text{for }t\in\left(  0,1\right)
\label{eigen1}\\
u^{\left(  p+j\right)  }\left(  0\right)   &  =u^{\left(  p+j\right)  }\left(
1\right)  =0\qquad\text{for }j=0,1,...,p-1. \label{eigen2}%
\end{align}
Kolmogorov used the following properties of problem (\ref{eigen1}%
)-(\ref{eigen2}) (cf. \cite{lorentz}, Chapter $9.6,$ Theorem $9,$ p. $146,$ or
\cite{naimark}, \cite{pinkus}):

\begin{proposition}
\label{PBVP}Problem (\ref{eigen1})-(\ref{eigen2}) has a countable set of
non-negative real eigenvalues with finite multiplicity. If we denote them by
$\lambda_{j}$ in a monotone order, they satisfy $\lambda_{j}\longrightarrow
\infty$ for $j\longrightarrow\infty.$ They satisfy the following asymptotic
$\lambda_{j}=\pi^{2p}j^{2p}\left(  1+O\left(  j^{-1}\right)  \right)  .$ The
corresponding orthonormalized eigenfunctions $\left\{  \psi_{j}\right\}
_{j=1}^{\infty}$ form a complete orthonormal system in $L_{2}\left(  \left[
0,1\right]  \right)  .$ The eigenvalue $\lambda=0$ has multiplicity $p$ and
the corresponding eigenfunctions $\left\{  \psi_{j}\right\}  _{j=1}^{p}$ are
the basis for the solutions to equation $u^{\left(  p\right)  }\left(
t\right)  =0$ in the interval $\left(  0,1\right)  .$
\end{proposition}

Further, Kolmogorov provided a description of the axes of the "cylindrical
ellipsoid" $K_{p}$, from which an approximation theorem of \emph{Jackson type
}easily follows (cf. \cite{lorentz}, chapter $4$ and chapter $5$).

\begin{proposition}
\label{PKolmogorovJackson} Let $f\in L_{2}\left(  \left[  a,b\right]  \right)
$ have the $L_{2}-$expansion
\[
f\left(  t\right)  =%
{\displaystyle\sum_{j=1}^{\infty}}
f_{j}\psi_{j}\left(  t\right)  .
\]
Then $f\in K_{p}$ if and only if
\[%
{\displaystyle\sum_{j=1}^{\infty}}
f_{j}^{2}\lambda_{j}\leq1.
\]
For $N\geq p+1$ and every $f\in K_{p}$ holds the following estimate
(\emph{Jackson type} approximation):
\begin{equation}
\left\Vert f-%
{\displaystyle\sum_{j=1}^{N}}
f_{j}\psi_{j}\left(  t\right)  \right\Vert _{L_{2}}\leq\frac{1}{\sqrt
{\lambda_{N+1}}}=O\left(  \frac{1}{\left(  N+1\right)  ^{p}}\right)  .
\label{Jackson1dim}%
\end{equation}

\end{proposition}

However, Kolmogorov didn't stop at this point but asked further, whether the
linear space $\widetilde{X}_{N}:=\left\{  \psi_{j}\right\}  _{j=1}^{N}$
provides the "best possible approximation among the linear spaces of dimension
$N$" in the following sense: If we put
\begin{equation}
d_{N}\left(  K_{p}\right)  :=\inf_{X_{N}}\operatorname*{dist}\left(
X_{N},K_{p}\right)  \label{KolmogorovWidth}%
\end{equation}
the main result he proved in \cite{kolmogorov1936} says
\begin{equation}
d_{N}\left(  K_{p}\right)  =\operatorname*{dist}\left(  \widetilde{X}%
_{N},K_{p}\right)  . \label{KolmogorovWidth2}%
\end{equation}
Here we have used the notations, to be used also further,
\begin{align}
\operatorname*{dist}\left(  X,K_{p}\right)   &  :=\sup_{y\in K_{p}%
}\operatorname*{dist}\left(  X,y\right) \label{dist1}\\
\operatorname*{dist}\left(  X,y\right)   &  =\inf_{x\in X}\left\Vert
x-y\right\Vert . \label{dist2}%
\end{align}
Hence, by inequality (\ref{Jackson1dim}), equality (\ref{KolmogorovWidth2})
reads as
\begin{align*}
d_{N}\left(  K_{p}\right)   &  =\frac{1}{\sqrt{\lambda_{N+1}}}\qquad\text{for
}N\geq p\\
d_{N}\left(  K_{p}\right)   &  =\infty\qquad\qquad\text{for }N=0,1,...,p-1.
\end{align*}

\begin{definition}
The left quantity in (\ref{KolmogorovWidth}) is called \textbf{Kolmogorov
}$N-$\textbf{width}, while the best approximation space $\widetilde{X}_{N}$ is
called \textbf{extremal (optimal) subspace} (cf. this terminology in
\cite{tikhomirov}, \cite{lorentz}, \cite{pinkus}).
\end{definition}

Thus the \emph{main approach to the successful application of the theory of
widths} is based on a Jackson type theorem by which a special space
$\widetilde{X}_{N}$ is identified. Then one has to find, among which subspaces
$X_{N}$ is $\widetilde{X}_{N}$ the extremal subspace. Put in a different
perspective : one has to find as wide class of spaces $X_{N}$ as possible,
among which $\widetilde{X}_{N}$ is the extremal subspace.

Now let us consider the following set which is a \emph{natural multivariate
generalization} of the above set $K_{p}$ defined in (\ref{Kp}): For a bounded
domain $B$ in $\mathbb{R}^{n}$ we put (more generally than (\ref{Kpstar}))
\begin{equation}
K_{p}^{\ast}:=\left\{  u\in H^{2p}\left(  B\right)  :%
{\displaystyle\int_{B}}
\left\vert L_{2p}u\left(  x\right)  \right\vert ^{2}dx\leq1\right\}  ,
\label{KpstarGeneral}%
\end{equation}
where $L_{2p}$ is a \emph{strongly elliptic} operator in $B.$ Let us remark
that the \emph{Sobolev space} $H^{2p}\left(  B\right)  $ is the multivariate
version of the space of absolutely continuous functions on the interval with a
highest derivative in $L_{2}$ (as in (\ref{Kp})). An important feature of the
set $K_{p}^{\ast}$ is that it contains an infinite-dimensional subspace
\[
\left\{  u\in H^{2p}\left(  B\right)  :L_{2p}u\left(  x\right)  =0,\quad
\text{for }x\in B\right\}  .
\]
Hence, all Kolmogorov widths are equal to infinity, i.e.
\begin{equation}
d_{N}\left(  K_{p}^{\ast}\right)  =\infty\qquad\text{for }N\geq0
\label{dNisInfinity}%
\end{equation}
and \emph{no way} is seen to improve this if one remains within the
finite-dimensional setting.

The main purpose of the present paper is to find a proper setting in the
framework of the Polyharmonic Paradigm which generalizes the above results of Kolmogorov.

\section{A reminder on Elliptic Boundary Value Problems \label{SreminderBVP}}

Let us specify the properties of the domains and the elliptic operators which
we will consider. In what follows we assume that the domain $D,$ the
differential operators and the boundary operators satisfy conditions for
\textbf{regular Elliptic BVP}. Namely, we give the following:

\begin{definition}
\label{Delliptic}We will say that the system of operators $\left\{
A;B_{j},\ j=1,2,...,m\right\}  $ forms a \textbf{regular Elliptic BVP in the
domain }$D\subset\mathbb{R}^{n}$ if the following conditions hold:

1. The operator
\[
A\left(  x,D_{x}\right)  =%
{\displaystyle\sum_{\left\vert \alpha\right\vert ,\left\vert \beta\right\vert
\leq m}}
\left(  -1\right)  ^{\left\vert \alpha\right\vert }D^{\alpha}a_{\alpha\beta
}\left(  x\right)  D^{\beta}%
\]
is a differential operator with a principal part defined as
\[
A_{0}\left(  x,D_{x}\right)  =%
{\displaystyle\sum_{\left\vert \alpha\right\vert +\left\vert \beta\right\vert
=2m}}
\left(  -1\right)  ^{\left\vert \alpha\right\vert }a_{\alpha\beta}\left(
x\right)  D^{\alpha+\beta}.
\]
It is \textbf{uniformly} \textbf{strongly elliptic}, i.e. for every $x\in D$
holds
\[
c_{0}\left\vert \xi\right\vert ^{2m}\leq\left\vert A_{0}\left(  x,\xi\right)
\right\vert \leq c_{1}\left\vert \xi\right\vert ^{2m}\qquad\text{for all real
}\xi\in\mathbb{R}^{n}\setminus\left\{  0\right\}  .
\]

2. The domain $D$ is bounded and has a boundary $\partial D$ of the class
$C^{2m}.$

3. For every pair of linearly independent real vectors $\xi,$ $\eta$ and
$x\in\overline{D}$ the polynomial in $z,$ $A_{0}\left(  x,\xi+z\eta\right)  $
has exactly $m$ roots with positive imaginary parts.

4. The coefficients of $A$ are in $C^{\infty}\left(  \overline{D}\right)  .$
The boundary operators $B_{j}\left(  x,D\right)  =%
{\displaystyle\sum_{\left\vert \alpha\right\vert \leq m_{j}}}
b_{j,\alpha}\left(  x\right)  D^{\alpha}$ form a \textbf{normal system}, i.e.
their principal symbols are \textbf{non-characteristic}, i.e. satisfy
$B_{j,0}\left(  x,\xi\right)  =%
{\displaystyle\sum_{\left\vert \alpha\right\vert =m_{j}}}
b_{j,\alpha}\left(  x\right)  \xi^{\alpha}\neq0$ for every $x\in\partial D$
and $\xi\neq0,$ $\xi$ is normal to $\partial D$ at $x;$ they have pairwise
different orders $m_{j}$ which satisfy $m_{j}<2m$ for $1\leq j\leq m,$ and
their coefficients $b_{j,\alpha}$ belong to $C^{\infty}$ in $\partial D.$

5. At any point $x\in\partial D$ let $\nu$ denote the outward normal to
$\partial D$ at $x$ and let $\xi\neq0$ be a real vector in the tangent
hyperplane to $\partial D$ at $x.$ The polynomials in $z$ given by
$B_{j,0}\left(  x,\xi+z\nu\right)  $ are linearly independent modulo the
polynomial $\prod_{k=1}^{m}\left(  z-z_{k}^{+}\left(  \xi\right)  \right)  $
where $z_{k}^{+}\left(  \xi\right)  $ denote the roots of $A_{0}\left(
x,\xi+z\eta\right)  $ with positive imaginary parts.
\end{definition}

\begin{remark}
With minor differences the above definition is available in
\cite{lions-magenes} (conditions (i)-(iii) in chapter $2,$ section $5.1$); in
\cite{taylor} (sections $5.11$ and $5.12$); in \cite{hoermanderVol3} (chapter
$20$); in \cite{okbook} (section $23.2,$ p. $473$).
\end{remark}

Let us define a special system of boundary operators called \emph{Dirichlet}.
We put
\begin{align*}
B_{j}  &  =\left(  \frac{\partial}{\partial n}\right)  ^{j-1}\qquad\text{for
}j=1,2,...,p-1\\
S_{j}  &  =\left(  \frac{\partial}{\partial n}\right)  ^{p+j-1}\qquad\text{for
}j=1,2,...,p-1.
\end{align*}
Obviously,
\[
\operatorname*{ord}\left(  B_{j}\right)  =j-1,\qquad\operatorname*{ord}\left(
S_{j}\right)  =p+j-1.
\]
Let us denote by $L_{2p}^{\ast}$ the operator formally adjoint to the elliptic
operator $L_{2p}$. There exist boundary operators $C_{j},$ $T_{j},$ for
$j=1,2,...,p-1,$ such that
\[
\operatorname*{ord}\left(  T_{j}\right)  =2p-j,\qquad\operatorname*{ord}%
\left(  C_{j}\right)  =p-j
\]
and the following Green's formula holds:\
\begin{equation}%
{\displaystyle\int_{B}}
\left(  L_{2p}u\cdot v-u\cdot L_{2p}^{\ast}v\right)  dx=%
{\displaystyle\sum_{j=0}^{p-1}}
{\displaystyle\int_{\partial B}}
\left(  S_{j}u\cdot C_{j}v-B_{j}u\cdot T_{j}v\right)  d\sigma_{y};
\label{GreenGeneral}%
\end{equation}
here $\partial_{n}$ denotes the normal derivative to $\partial B,$ for
functions $u$ and $v$ in the classes of Sobolev, $u,v$ $\in$ $H^{2p}\left(
B\right)  $ (cf. \cite{lions-magenes}, Theorem $2.1$ in section $2.2,$ chapter
$2$, and Remark $2.2$ in section $2.3$).

For us the following eigenvalue problem will be important to consider for
$U\in H^{2p}\left(  B\right)  ,$ which is analogous to problem (\ref{eigen1}%
)-(\ref{eigen2}):%
\begin{align}
L_{2p}^{\ast}L_{2p}U\left(  x\right)   &  =\lambda U\left(  x\right)
\qquad\qquad\qquad\text{for }x\in B\label{eigen1Multi}\\
B_{j}L_{2p}U\left(  y\right)   &  =S_{j}L_{2p}U\left(  y\right)
=0,\qquad\text{for }y\in\partial B,\quad j=0,1,...,p-1 \label{eigen2Multi}%
\end{align}
where $\partial_{n}$ denotes the normal derivative at $y\in\partial B.$ It is
obvious that the operator $L_{2p}^{\ast}L_{2p}$ is formally self-adjoint,
however the BVP (\ref{eigen1Multi})-(\ref{eigen2Multi}) is not a nice one.
Since a direct reference seems not to be available, we provide its
consideration in the following theorem which is an analog to Proposition
\ref{PBVP}.

\begin{theorem}
\label{TExpansionBerezanskii} Let the operator $L_{2p}$ be uniformly strongly
elliptic in the domain $B.$ Then problem (\ref{eigen1Multi}%
)-(\ref{eigen2Multi}) has only real non-negative eigenvalues.

1. The eigenvalue $\lambda=0$ has infinite multiplicity with corresponding
eigenfunctions $\left\{  \psi_{j}^{\prime}\right\}  _{j=1}^{\infty}$ which
represent an orthonormal basis of the space of all solutions to the equation
$L_{2p}U\left(  x\right)  =0,$ for $x\in B.$

2. The positive eigenvalues are countably many and each has \textbf{finite
multiplicity}, and if we denote them by $\lambda_{j}$ ordered increasingly,
they satisfy $\lambda_{j}\longrightarrow\infty$ for $j\longrightarrow\infty.$

3. The orthonormalized eigenfunctions, corresponding to eigenvalues
$\lambda_{j}>0,$ will be denoted by $\left\{  \psi_{j}\right\}  _{j=1}%
^{\infty}.$ The set of functions $\left\{  \psi_{j}\right\}  _{j=1}^{\infty
}\bigcup\left\{  \psi_{j}^{\prime}\right\}  _{j=1}^{\infty}$ form a complete
orthonormal system in $L_{2}\left(  B\right)  .$
\end{theorem}

\begin{remark}
Problem (\ref{eigen1Multi})-(\ref{eigen2Multi}) is well known to be a
non-regular elliptic BVP, as well as non-coercive variational, cf.
\cite{agmon} ( p. $150$ ) and \cite{lions-magenes} (Remark $9.8$ in chapter
$2,$ section $9.6$, and section $9.8$ ).
\end{remark}

The proof is provided in the Appendix below, section \ref{Sappendix}.

\section{The principal axes of the ellipsoid $K_{p}^{\ast}$ and a Jackson type
theorem \label{SprincipalAxes}}

Here we will find the principal exes of the ellipsoid $K_{p}^{\ast}$ defined
as
\begin{equation}
K_{p}^{\ast}:=\left\{  u\in H^{2p}\left(  B\right)  :%
{\displaystyle\int_{B}}
\left\vert L_{2p}u\left(  x\right)  \right\vert ^{2}dx\leq1\right\}  ,
\label{KpstarGENERAL}%
\end{equation}
where $L_{2p}$ is a \emph{uniformly strongly elliptic} operator in $B.$

We prove the following theorem which generalizes Kolmogorov's one-dimensional
result from Proposition \ref{PKolmogorovJackson}, about the representation of
the ellipsoid $K_{p}$ in principal axes.

\begin{theorem}
\label{TprincipalAxesMultivariate}Let $f\in K_{p}^{\ast}.$ Then $f$ is
represented in a $L_{2}-$series as
\[
f\left(  x\right)  =%
{\displaystyle\sum_{j=1}^{\infty}}
f_{j}^{\prime}\psi_{j}^{\prime}\left(  x\right)  +%
{\displaystyle\sum_{j=1}^{\infty}}
f_{j}\psi_{j}\left(  x\right)  ,
\]
where by Theorem \ref{TExpansionBerezanskii} the eigenfunctions $\psi
_{j}^{\prime}$ satisfy $\Delta^{p}\psi_{j}^{\prime}\left(  x\right)  =0$ while
the eigenfunctions $\psi_{j}$ correspond to the eigenvalues $\lambda_{j}>0,$
and also
\begin{equation}%
{\displaystyle\sum_{j=1}^{\infty}}
\lambda_{j}f_{j}^{2}\leq1. \label{EllipsoidCondition}%
\end{equation}
Vice versa, every sequence $\left\{  f_{j}^{\prime}\right\}  _{j=1}^{\infty
}\bigcup\left\{  f_{j}\right\}  _{j=1}^{\infty}$ with $%
{\displaystyle\sum_{j=1}^{\infty}}
\left\vert f_{j}^{\prime}\right\vert ^{2}+%
{\displaystyle\sum_{j=1}^{\infty}}
\left\vert f_{j}\right\vert ^{2}<\infty$ and $%
{\displaystyle\sum_{j=1}^{\infty}}
\lambda_{j}f_{j}^{2}\leq1$ defines a function $f\in L_{2}\left(  B\right)  $
which is in $K_{p}^{\ast}.$
\end{theorem}

%

\proof
\textbf{(1)} According to Theorem \ref{TExpansionBerezanskii}, we know that
arbitrary $f\in L_{2}\left(  B\right)  $ is represented as
\begin{align*}
f\left(  x\right)   &  =%
{\displaystyle\sum_{j=1}^{\infty}}
f_{j}^{\prime}\psi_{j}^{\prime}\left(  x\right)  +%
{\displaystyle\sum_{j=1}^{\infty}}
f_{j}\psi_{j}\left(  x\right) \\
\left\Vert f\right\Vert _{L_{2}}^{2}  &  =%
{\displaystyle\sum_{j=1}^{\infty}}
\left\vert f_{j}^{\prime}\right\vert ^{2}+%
{\displaystyle\sum_{j=1}^{\infty}}
\left\vert f_{j}\right\vert ^{2}<\infty
\end{align*}
with convergence in the space $L_{2}\left(  B\right)  .$

\textbf{(2)} From the proof of Theorem \ref{TExpansionBerezanskii}, we know
that if we put
\[
\phi_{j}\left(  x\right)  =L_{2p}\psi_{j}\left(  x\right)  \qquad\text{for
}j\geq1,
\]
then the system of functions
\[
\frac{\phi_{j}\left(  x\right)  }{\sqrt{\lambda_{j}}}\qquad\text{for }j\geq1
\]
is orthonormal sequence which is complete in $L_{2}\left(  B\right)  .$

\textbf{(3)} We will prove now that if $f\in L_{2}\left(  B\right)  $ then
$f\in K_{p}^{\ast}$ iff
\[%
{\displaystyle\sum_{j=1}^{\infty}}
f_{j}^{2}\lambda_{j}\leq1.
\]
Indeed, for every $f\in H^{2p}\left(  B\right)  $ we have the expansion
$f\left(  x\right)  =%
{\displaystyle\sum_{j=1}^{\infty}}
f_{j}^{\prime}\psi_{j}^{\prime}\left(  x\right)  +%
{\displaystyle\sum_{j=1}^{\infty}}
f_{j}\psi_{j}\left(  x\right)  .$ We want to see that it is possible to
differentiate termwise this expansion, i.e.
\[
L_{2p}f\left(  x\right)  =%
{\displaystyle\sum_{j=1}^{\infty}}
f_{j}L_{2p}\psi_{j}\left(  x\right)  =%
{\displaystyle\sum_{j=1}^{\infty}}
f_{j}\phi_{j}\left(  x\right)
\]
Since $\left\{  \frac{\phi_{j}}{\sqrt{\lambda_{j}}}\right\}  _{j\geq1}$ is a
complete orthonormal basis of $L_{2}\left(  B\right)  $ it is sufficient to
see that
\[%
{\displaystyle\int_{B}}
L_{2p}f\left(  x\right)  \phi_{j}dx=%
{\displaystyle\int_{B}}
\left(
{\displaystyle\sum_{j=1}^{\infty}}
f_{j}L_{2p}\psi_{j}\left(  x\right)  \right)  \phi_{j}dx.
\]
Due to the boundary properties of $\phi_{j}$ and since $\phi_{j}=L_{2p}%
\psi_{j},$ we obtain
\[%
{\displaystyle\int_{B}}
L_{2p}f\left(  x\right)  \phi_{j}dx=%
{\displaystyle\int_{B}}
f\left(  x\right)  L_{2p}^{\ast}\phi_{j}dx=\lambda_{j}%
{\displaystyle\int_{B}}
f\psi_{j}dx=\lambda_{j}f_{j}.
\]
On the other hand
\[%
{\displaystyle\int_{B}}
\left(
{\displaystyle\sum_{k=1}^{\infty}}
f_{k}\phi_{k}\left(  x\right)  \right)  \phi_{j}dx=\lambda_{j}f_{j}.
\]
Hence
\[
L_{2p}f\left(  x\right)  =%
{\displaystyle\sum_{j=1}^{\infty}}
f_{j}L_{2p}\psi_{j}\left(  x\right)  =%
{\displaystyle\sum_{j=1}^{\infty}}
f_{j}\phi_{j}\left(  x\right)  =%
{\displaystyle\sum_{j=1}^{\infty}}
\sqrt{\lambda_{j}}f_{j}\frac{\phi_{j}\left(  x\right)  }{\sqrt{\lambda_{j}}}%
\]
and since $\left\{  \frac{\phi_{j}}{\sqrt{\lambda_{j}}}\right\}  _{j\geq1}$ is
an orthonormal system, it follows
\[
\left\Vert L_{2p}f\right\Vert _{L_{2}}^{2}=%
{\displaystyle\sum_{j=1}^{\infty}}
\lambda_{j}f_{j}^{2}.
\]
Thus if $f\in K_{p}$ it follows that $%
{\displaystyle\sum_{j=1}^{\infty}}
\lambda_{j}f_{j}^{2}\leq1.$

Now, assume vice versa, that $%
{\displaystyle\sum_{j=1}^{\infty}}
f_{j}^{2}\lambda_{j}\leq1$ holds together with $%
{\displaystyle\sum_{j=1}^{\infty}}
\left\vert f_{j}^{\prime}\right\vert ^{2}+%
{\displaystyle\sum_{j=1}^{\infty}}
\left\vert f_{j}\right\vert ^{2}<\infty$. We have to see that the function
\[
f\left(  x\right)  =%
{\displaystyle\sum_{j=1}^{\infty}}
f_{j}^{\prime}\psi_{j}^{\prime}\left(  x\right)  +%
{\displaystyle\sum_{j=1}^{\infty}}
f_{j}\psi_{j}\left(  x\right)
\]
belongs to the space $H^{2p}\left(  B\right)  .$ Based on the completeness and
orthonormality of the system $\left\{  \frac{\phi_{j}\left(  x\right)  }%
{\sqrt{\lambda_{j}}}\right\}  _{j=1}^{\infty}$ we may define the function
$g\in L_{2}$ by putting
\[
g\left(  x\right)  =%
{\displaystyle\sum_{j=1}^{\infty}}
\sqrt{\lambda_{j}}f_{j}\frac{\phi_{j}\left(  x\right)  }{\sqrt{\lambda_{j}}}=%
{\displaystyle\sum_{j=1}^{\infty}}
f_{j}\phi_{j}\left(  x\right)  ;
\]
it obviously satisfies $\left\Vert g\right\Vert _{L_{2}}\leq1.$

From the local solvability of elliptic equations (\cite{lions-magenes}) there
exists a function $F\in H^{2p}\left(  B\right)  $ which is a solution to
equation $L_{2p}F=g.$ Let its representation be
\[
F\left(  x\right)  =%
{\displaystyle\sum_{j=1}^{\infty}}
f_{j}^{\prime}\psi_{j}^{\prime}\left(  x\right)  +%
{\displaystyle\sum_{j=1}^{\infty}}
F_{j}\psi_{j}\left(  x\right)
\]
with some coefficients $F_{j}$ satisfying $%
{\displaystyle\sum_{j}}
\left\vert F_{j}\right\vert ^{2}<\infty.$ As above we obtain
\begin{align*}
\lambda_{j}%
{\displaystyle\int_{B}}
F\psi_{j}dx  &  =%
{\displaystyle\int_{B}}
FL_{2p}^{\ast}L_{2p}\psi_{j}dx=%
{\displaystyle\int_{B}}
L_{2p}F\cdot L_{2p}\psi_{j}dx\\
&  =%
{\displaystyle\int_{B}}
g\cdot\phi_{j}dx
\end{align*}
which implies $F_{j}=f_{j}.$ Hence, $F=f$ and $f\in H^{2p}\left(  B\right)  .$
This ends the proof.%

\endproof

We are able to prove finally a \emph{Jackson type} result as in Proposition
\ref{PKolmogorovJackson}.

\begin{theorem}
Let $N\geq1.$ Then for every $N\geq1$ and every $f\in K_{p}^{\ast}$ holds the
following estimate:
\[
\left\Vert f-%
{\displaystyle\sum_{j=1}^{\infty}}
f_{j}^{\prime}\psi_{j}^{\prime}\left(  x\right)  -%
{\displaystyle\sum_{j=1}^{N}}
f_{j}\psi_{j}\left(  x\right)  \right\Vert _{L_{2}}\leq\frac{1}{\sqrt
{\lambda_{N+1}}}.
\]

\end{theorem}

%

\proof
The proof follows directly. Indeed, due to the monotonicity of $\lambda_{j},$
and inequality (\ref{EllipsoidCondition}), we obtain
\[
\left\Vert f-%
{\displaystyle\sum_{j=1}^{\infty}}
f_{j}^{\prime}\psi_{j}^{\prime}\left(  x\right)  -%
{\displaystyle\sum_{j=1}^{N}}
f_{j}\psi_{j}\left(  x\right)  \right\Vert _{L_{2}}^{2}=%
{\displaystyle\sum_{j=N+1}^{\infty}}
f_{j}^{2}\leq\frac{1}{\lambda_{N+1}}%
{\displaystyle\sum_{j=N+1}^{\infty}}
f_{j}^{2}\lambda_{j}\leq\frac{1}{\lambda_{N+1}}.
\]
This ends the proof.%

\endproof

\section{Introducing the Hierarchy and Harmonic Widths \label{Shierarchy}}

In the present section we introduce the simplest representatives of the class
of domains having Harmonic Dimension $N,$ which are called \textbf{First Kind}
domains. They are piece-wise solutions to regular elliptic equations$.$

\begin{definition}
\label{Dhdimension} Let $D\subset\mathbb{R}^{n}$ be a bounded domain. For an
integer $M\geq1$ we say that the linear subspace $X_{M}\subset L_{2}\left(
D\right)  $ is of \textbf{First Kind} and has \textbf{Harmonic Dimension} $M,$
and write
\begin{equation}
\operatorname{hdim}\left(  X_{M}\right)  =M, \label{hdim}%
\end{equation}
if the following conditions are fulfilled:

1. There exists a finite number of domains $D_{j}$ with \textbf{piece-wise
smooth} boundaries $\partial D_{j}$ (which guarantees the validity of Green's
formula (\ref{GreenGeneral})), which are pairwise disjoint, i.e. $D_{i}\bigcap
D_{j}=\varnothing$ for $i\neq j,$ and such that we have the \textbf{domain
partition}
\begin{equation}
D=\bigcup_{j}D_{j}. \label{domainPartition}%
\end{equation}

2. We assume that for $j=1,2,...,k$ the \textbf{factorization operators}
$Q_{j}\left(  x,D_{x}\right)  $ are \textbf{uniformly strongly elliptic} in
the domain $D$ and, the functions $\rho_{j}$ defined in $D$ are infinitely
smooth and satisfy
\[
Z:=\bigcup_{j=1}^{k}\left\{  x\in D:\rho_{j}\left(  x\right)  =0\right\}
\subset\bigcup_{j=1}^{k}\partial D_{j}\setminus\partial D.
\]
We assume that $ord\left(  Q_{j}\right)  =2N_{j}$ and $%
{\displaystyle\sum_{j=1}^{k}}
2N_{j}=2M.$ Define the operator
\begin{equation}
P_{2M}\left(  x,D_{x}\right)  u\left(  x\right)  =\left(  \prod_{j=1}^{k}%
Q_{j}\left(  x,D_{x}\right)  \frac{1}{\rho_{j}\left(  x\right)  }\right)
u\left(  x\right)  \label{P2Nproduct}%
\end{equation}
for the points $x\in D$ where it is correctly defined (out of the set $Z$ ).

We specify the interface conditions:\ Let us denote by $u_{i}=u_{|D_{i}}$ the
restriction of $u$ to $D_{i}.$ If for some indexes $i\neq j,\ $the
intersection $H:=\partial D_{i}\bigcap\partial D_{j}$ has nonempty interior in
the relative topology of $\partial D_{i}$ (hence also in $\partial D_{j}$)
then the following \textbf{interface conditions} hold on $H$ in the sense of
traces:
\begin{equation}
\left(  \frac{\partial}{\partial n_{x}}\right)  ^{k}u_{i}\left(  x\right)
=\left(  \frac{\partial}{\partial n_{x}}\right)  ^{k}u_{j}\left(  x\right)
\qquad\text{for }k=0,1,...,2M-1; \label{Interfaces}%
\end{equation}
here the vector $n_{x}$ denotes one of the normals at $x$ to the surface
$\partial D_{i}\bigcap\partial D_{j}.$

We define the space $X_{M}$ by putting
\begin{equation}
X_{M}=\left\{
\begin{array}
[c]{c}%
u\in H^{2M}\left(  D\right)  :P_{2M}u\left(  x\right)  =0,\quad\text{for }%
x\in\bigcup_{j=1}^{k}D_{j},\ \\
\text{and }u\text{ satisfies the interface conditions (\ref{Interfaces})}%
\end{array}
\right\}  . \label{XNDj}%
\end{equation}

\end{definition}

\begin{remark}
1. In \cite{kounchevSozopol} we considered the case of spaces $X_{N}$ of
Harmonic Dimension $N$ defined by a single elliptic operator $P_{2N}$ (i.e.
$P_{2N}=Q_{1}$) and a trivial partition of $D,$ i.e. $D=D_{1}.$

2. Let us comment on the interface conditions (\ref{Interfaces}) in Definition
\ref{Dhdimension}. Let us assume that we have an elliptic operator $P_{2N}$
with smooth coefficients defined on $D$ and that a non-trivial partition
$\bigcup D_{j}$ is given. Due to the piece-wise smoothness of the boundaries
$\partial D_{j}$ we may apply the Green formula, and from the interface
conditions (\ref{Interfaces}) it follows that "analytic continuation" is
possible, hence every function in $X_{N}$ is a solution to $P_{2N}u=0$ in the
whole domain $D$ (see similar result in \cite{okbook}, Lemma $20.10,$ and the
proof of Theorem $20.11$).

3. One may choose a different set of interface conditions which are equivalent
to (\ref{Interfaces}), see \cite{okbook} (Remark $20.12$), and
\cite{lions-magenes} (Lemma $2.1$ in chapter $2$).

4. The spaces $X_{N}$ defined in Definition \ref{Dhdimension} mimic in a
natural way the one-dimensional case: the operator $P_{2M}$ (\ref{P2Nproduct})
is similar to the operator (\ref{LNoperatorForm}) in Proposition
\ref{PstructureMarkovSystem}.

5. The operator $\prod_{j}\rho_{j}\cdot P_{2M}$ does not have a singularity in
the principal symbol but eventually only in the lower order coefficients.
\end{remark}

Here is a simple non-trivial example to Definition \ref{Dhdimension}:
\begin{align*}
D_{1}  &  =\left\{  x:\left\vert x\right\vert <1\right\}  ,\quad
D_{2}=\left\{  x:1<\left\vert x\right\vert <2\right\} \\
D  &  =\left\{  x:\left\vert x\right\vert <2\right\} \\
P_{4}^{1}\left(  x;D_{x}\right)  u\left(  x\right)   &  =\Delta\frac
{1}{1-\left\vert x\right\vert }\Delta u\left(  x\right)  \qquad\text{for }x\in
D_{1}\\
P_{4}^{2}\left(  x;D_{x}\right)  u\left(  x\right)   &  =-\Delta\frac
{1}{1-\left\vert x\right\vert }\Delta u\left(  x\right)  \qquad\text{for }x\in
D_{2},
\end{align*}
where $\Delta$ is the Laplace operator. Typical elements of $X_{2}$ are the
functions $u$ which are obtained as solutions to
\[
\Delta u=\left(  1-\left\vert x\right\vert \right)  w\qquad\text{in }D,
\]
where $\Delta w=0$ in $D.$

The following result shows that we may construct a lot of solutions belonging
to the set $X_{M}$ of Definition \ref{Dhdimension}. We call these "direct solutions".

\begin{proposition}
\label{Pdirectsolutions}There is a set of boundary conditions $B_{\ell},$
$\ell=1,2,...,M$ on $\partial D$ such that problem
\begin{align*}
P_{2M}u\left(  x\right)   &  =0\qquad\text{for }x\in D\\
B_{\ell}u\left(  y\right)   &  =h_{\ell}\left(  y\right)  \qquad\text{for
}y\in\partial D,\ \text{and }\ell=1,2,...,M
\end{align*}
is solvable for arbitrary data $\left\{  h_{\ell}\right\}  _{\ell=1}^{M}$ from
the corresponding Sobolev spaces, i.e. $h_{\ell}\in H^{2M-ord\left(  B_{\ell
}\right)  -1/2}\left(  \partial D\right)  ,$ and the solution has the maximal
regularity, i.e. $u\in H^{2M}\left(  D\right)  .$
\end{proposition}

%

\proof
For every $j$ with $1\leq j\leq k$ we choose boundary operators $B_{j,m}$ for
$m=1,2,...,N_{j}$ for a \emph{regular elliptic BVP} $\left\{  Q_{j}\left(
x,D_{x}\right)  ;B_{j,m},\ m=1,2,...,N_{j}\right\}  .$ If we have the data
function $f$ on $D$ and $h^{\left(  j\right)  }=\left\{  h_{j,m}\right\}
_{m=1}^{N_{j}}$ on the boundary $\partial D$ then the solution of the elliptic
BVP
\begin{align}
Q_{j}w  &  =f\qquad\text{on }D\label{Qjw}\\
B_{j,m}w  &  =h_{j,m},\qquad\text{on }\partial D,\ \text{for }m=1,2,...,N_{j}
\label{Qjw2}%
\end{align}
in case it exists will be denoted by $I_{j}\left(  f,h^{\left(  j\right)
}\right)  .$\footnote{For the solvability recall that there is a finite number
of conditions which have to be satisfied by the data $\left\{  f,h_{j,m}%
\right\}  $ which guarantee the solvability, cf. \cite{lions-magenes} (Theorem
$5.3,$ chapter $2,$ section $5.3$).} We may write inductively
\[
u=\rho_{k}I_{k}\left(  \cdot\cdot\cdot\rho_{2}I_{2}\left(  \rho_{1}%
I_{1}\left(  0;h^{\left(  1\right)  }\right)  ;h^{\left(  2\right)  }\right)
\cdot\cdot\cdot\right)  .
\]

For simplicity of notation let us assume that $k=2.$ Then the boundary
conditions satisfied by $u$ are obtained from
\begin{align*}
Q_{1}w  &  =0\\
B_{1,m}w  &  =h_{1,m}\qquad\text{for }m=1,2,...,N_{1}%
\end{align*}
and
\begin{align*}
Q_{2}\left(  \frac{1}{\rho_{2}}u\right)   &  =\rho_{1}w\\
B_{2,m}\frac{1}{\rho_{2}}u  &  =h_{2,m}\qquad\text{for }m=1,2,...,N_{2}%
\end{align*}
hence, we obtain
\[
B_{1,m}\left(  \frac{1}{\rho_{1}}Q_{2}\left(  \frac{1}{\rho_{2}}u\right)
\right)  =h_{1,m}\qquad\text{for }m=1,2,...,N_{1}.
\]
Thus we see that the system of boundary operators on $\partial D$
\begin{align*}
B_{2,m}\frac{1}{\rho_{2}}u,\qquad\text{for }m  &  =1,2,...,N_{2}\\
B_{1,m}\left(  \frac{1}{\rho_{1}}Q_{2}\left(  \frac{1}{\rho_{2}}u\right)
\right)  \qquad\text{for }m  &  =1,2,...,N_{1}%
\end{align*}
is normal. Let us put
\begin{align*}
B_{j}u  &  =B_{2,j}\frac{1}{\rho_{2}}u\qquad\qquad\qquad\qquad\text{for
}j=1,2,...,N_{2}\\
B_{N_{2}+j}u  &  =B_{1,j}\left(  \frac{1}{\rho_{1}}Q_{2}\left(  \frac{1}%
{\rho_{2}}u\right)  \right)  \qquad\text{for }j=1,2,...,N_{1}.
\end{align*}
A simple direct check shows that the orders of the system of operators
\[
\left\{  B_{j}:j=1,2,...,N_{1}+N_{2}\right\}
\]
differ, and also satisfy the condition for being "non-characteristic" on the
boundary, cf. Definition \ref{Delliptic}, item 4). We may proceed inductively
to prove the statement for arbitrary $k\geq3.$%

\endproof

\begin{remark}
Apparently, one may prove that the set of "direct solutions" obtained in
Proposition \ref{Pdirectsolutions} is dense in the whole space $X_{M}$ defined
in Definition \ref{Dhdimension}.
\end{remark}

The following fundamental theorem shows that, as in the one-dimensional case,
on arbitrary small sub-domain $G$ in $D$ with $G\bigcap\left(  \bigcup\partial
D_{j}\right)  =\varnothing,$ the space $X_{M}$ with $\operatorname{hdim}%
\left(  X_{M}\right)  =M$ has the same Harmonic Dimension $M.$ From a
different point of view, it shows that a theorem of Runge-Lax-Malgrange type
is true also for elliptic operators with singular coefficients of the type of
operators $P_{2M}$ considered in Definition \ref{Dhdimension}.

\begin{theorem}
\label{Tdensity} Let the First Kind space $X_{M}$ satisfy Definition
\ref{Dhdimension} with
\[
\operatorname{hdim}\left(  X_{M}\right)  =M.
\]
Assume that the elliptic operator $P_{2M}$ which corresponds to the space
$X_{M}$ has factorization operators $Q_{j}$ (from (\ref{P2Nproduct}))
satisfying condition $\left(  U\right)  _{s}$ for uniqueness in the Cauchy
problem in the small.\footnote{The differential operator $P$ satisfies
condition $\left(  U\right)  _{s}$ for uniqueness in the Cauchy problem in the
small in $G$ provided that if $G_{1}$ is a connected open subset of $G$ and
$u\in C^{r}\left(  G_{1}\right)  $ is a solution to $P^{\ast}u=0$ and $u$ is
zero on a non-emplty subset of $G_{1}$ then $u$ is identically zero. Elliptic
operators with analytic coefficients satisfy this property (cf.
\cite{bersJohnschechter}, part $II,$ chapter $1.4;$ \cite{browder}, p.
$402$).} Let $G$ be a compact subdomain in some $D_{j},$ i.e. $G\bigcap\left(
\bigcup\partial D_{j}\right)  =\varnothing.$ Then the set of "direct
solutions" considered in Proposition \ref{Pdirectsolutions} is dense in
$L_{2}\left(  G\right)  $ in the space
\[
\left\{  u\in H^{2M}\left(  G\right)  :P_{2M}u=0\quad\text{in }G\right\}  .
\]

\end{theorem}

%

\proof
For simplicity of notations we assume that for the elliptic operator $P_{2M}$
associated with $X_{M},$ by Definition \ref{Dhdimension}, we have only two
factorizing operators $Q_{1}$ and $Q_{2},$ i.e. $P_{2M}u=Q_{1}\frac{1}%
{\rho_{1}}Q_{2}\left(  \frac{1}{\rho_{2}}u\right)  .$

Let us take a solution $u\in H^{2M}\left(  G\right)  $ to $P_{2M}u=0$ in $G.$
We have
\[
Q_{1}\frac{1}{\rho_{1}}Q_{2}\left(  \frac{1}{\rho_{2}}u\right)  =0\qquad
\text{in }G
\]
and we use the solutions $I_{j}$ for the Elliptic BVP (\ref{Qjw})-(\ref{Qjw2})
considered in the domain $G,$ to express arbitrary solution as
\[
u=\rho_{2}I_{2}\left(  \rho_{1}I_{1}\left(  0;h^{\left(  1\right)  }\right)
;h^{\left(  2\right)  }\right)  ,
\]
where the boundary data $h^{\left(  1\right)  }$ and $h^{\left(  2\right)  }$
are arbitrary in proper Sobolev spaces. By the approximation theorem of
Runge-Lax-Malgrange type (cf. \cite{browder}, Theorem $4,$ and references
there), which uses essentially property $\left(  U\right)  _{s}$ of operator
$Q_{1}^{\ast},$ we obtain a function $w_{\varepsilon}$ which is a solution to
$Q_{1}w_{\varepsilon}=0$ in $D$ and such that
\[
\left\Vert I_{1}\left(  0;h^{\left(  1\right)  }\right)  -w_{\varepsilon
}\right\Vert _{L_{2}\left(  G\right)  }<\varepsilon.
\]
Next we apply the same approximation argument but with non-zero right-hand
side $\rho_{1}I_{1}\left(  0;h^{\left(  1\right)  }\right)  $ (cf.
\cite{browder2}) to prove the existence of a function $v_{\varepsilon}$ such
that
\[
\left\Vert I_{2}\left(  \rho_{1}I_{1}\left(  0;h^{\left(  1\right)  }\right)
;h^{\left(  2\right)  }\right)  -v_{\varepsilon}\right\Vert <C\varepsilon
\]
for some constant $C>0,$ where the constant $C$ depends on the functions
$\rho_{j}.$ Thus we obtain the function
\[
u_{\varepsilon}=\rho_{2}v_{\varepsilon}%
\]
which satisfies
\[
\left\Vert u_{\varepsilon}-u\right\Vert _{L_{2}\left(  G\right)  }%
<C_{1}\varepsilon,
\]
and is a "direct solution" in the sense of Proposition \ref{Pdirectsolutions}.%

\endproof

The following theorem studies the orthogonal complement $X_{N}\ominus X_{M}$
of two First Kind spaces where $M<N.$ While we will not need the whole
generality of the result proved, the proof shows that $X_{N}\ominus X_{M}$ has
at least $\operatorname{hdim}$ equal to $N-M.$

\begin{theorem}
\label{TtransversalSpaces}Let $M<N$ and the First Kind spaces $X_{M},$ $X_{N}$
satisfy Definition \ref{Dhdimension} with
\[
\operatorname{hdim}\left(  X_{M}\right)  =M,\quad\operatorname{hdim}\left(
X_{N}\right)  =N.
\]
Assume that the elliptic operator $P_{2N}^{\prime},$ which is associated with
the space $X_{N},$ has (by (\ref{P2Nproduct})) factorization operators $Q_{j}$
satisfying condition $\left(  U\right)  _{s}$ for uniqueness in the Cauchy
problem in the small (as in Theorem \ref{Tdensity}). Then the space
$Y=X_{N}\setminus X_{M}$ is infinite-dimensional.
\end{theorem}

%

\proof
\textbf{(1)} Let, by Definition \ref{Dhdimension}, the partition $\bigcup
D_{j}$ and the operator $P_{2M}$ correspond to $X_{M},$ while the partition
$\bigcup D_{j}^{\prime}$ and the operator $P_{2N}^{\prime}$ correspond to
$X_{N}.$ Assume that $D_{1}\bigcap D_{1}^{\prime}\neq\varnothing.$ Then we
will choose a subdomain $G$ which is compactly supported in $D_{1}\bigcap
D_{1}^{\prime}.$

Further we will fix our attention to the subdomain $G$ where both operators
$P_{2M}$ and $P_{2N}^{\prime}$ are uniformly strongly elliptic and will
construct a subset of $X_{N}\ominus X_{M}$ restricted to the domain $G.$ Let
us be more precise: If we denote by
\begin{equation}
X_{N}^{G}:=\left\{  u:H^{2N}\left(  G\right)  :P_{2N}^{\prime}u=0\quad\text{in
}G\right\}  \label{XNG}%
\end{equation}
then we will construct an infinite-dimensional subspace of $X_{N}%
^{G}\circleddash X_{M}^{G}.$

\textbf{(2)} For the \emph{uniformly strongly elliptic} operator $P_{2M}$ on
the domain $G$ we choose the Dirichlet system of boundary operators
$B_{j}=\frac{\partial^{j-1}}{\partial n^{j-1}},$ for $j\geq1,$ which are
iterates of the normal derivative $\frac{\partial}{\partial n}$ on the
boundary $\partial G.$ As already mentioned the system of operators $\left\{
P_{2M};\frac{\partial^{j}}{\partial n^{j}}:j=0,1,...,M-1\right\}  $ on $G$
forms a \emph{regular Elliptic BVP }(this is the Dirichlet Elliptic BVP for
the operator $P_{2M}$) (cf. \cite{lions-magenes}, Remark $1.3$ in section
$1.4$, chapter $2$).

We complete the system $\left\{  B_{j}\right\}  _{j=1}^{M}$ by the system of
boundary operators $S_{j}=\frac{\partial^{M-1+j}}{\partial n^{M-1+j}}$ for
$j=1,2,...M.$ Hence, the system composed $\left\{  B_{j}\right\}  _{j=1}%
^{M}\bigcup\left\{  S_{j}\right\}  _{j=1}^{M}$ is a \emph{Dirichlet system} of
order $2M$ (cf. \cite{lions-magenes}, Definition $2.1$ and Theorem $2.1$ in
section $2.2,$ chapter $2$). Further, by \cite{lions-magenes} (Theorem $2.1$),
there exists a unique Dirichlet system of order $2M$ of boundary operators
$\left\{  C_{j},T_{j}\right\}  _{j=1}^{M}$ which is uniquely determined as the
adjoint to the system $\left\{  B_{j},S_{j}\right\}  _{j=1}^{M},$ and the
Green formula (\ref{GreenGeneral}) holds on the domain $G.$ We will use this below.

\textbf{(3)} In the domain $G$ we consider the elliptic operator
$P_{2N}^{\prime}P_{2M}^{\ast}.$ As a product of two \emph{ strongly elliptic}
operators it is such again. By a standard construction cited above (cf.
\cite{lions-magenes}, Theorem $2.1,$ section $2.2,$ chapter $2$), we may
complete the Dirichlet system of operators $\left\{  B_{j},S_{j}\right\}
_{j=1}^{M}$ with $N-M$ boundary operators $R_{j}=\frac{\partial^{2M-1+j}%
}{\partial n^{2M-1+j}},$ $j=1,2,...,N-M.$ Again by the above cited theorem,
the Dirichlet system of boundary operators
\[
\left\{  B_{j},S_{j}\right\}  _{j=1}^{M}%
{\textstyle\bigcup}
\left\{  R_{j}\right\}  _{j=1}^{N-M}%
\]
\emph{covers} the operator $P_{2N}^{\prime}P_{2M}^{\ast}.$ Finally, we
consider the solutions $g\in H^{2N+2M}\left(  G\right)  $ to the following
Elliptic BVP:
\begin{align}
P_{2N}^{\prime}P_{2M}^{\ast}g\left(  x\right)   &  =0\qquad\qquad
\qquad\ \ \text{for }x\in G\label{gorthogonal1}\\
B_{j}g\left(  y\right)   &  =S_{j}g\left(  y\right)  =0\qquad\text{for
}j=0,1,...,N-1,\text{ for }y\in\partial G\label{gorthogonal2}\\
R_{j}g\left(  y\right)   &  =h_{j}\left(  y\right)  \qquad\qquad\ \ \text{for
}j=1,2,...,N-M,\text{ for }y\in\partial G. \label{gorthogonal3}%
\end{align}
We may apply a classical result \cite{lions-magenes} (the existence Theorem
$5.2$ and Theorem $5.3$ in chapter $2$), to the solvability of problem
(\ref{gorthogonal1})-(\ref{gorthogonal3}) in the space $H^{2M+2N}\left(
G\right)  .$

\textbf{(4) }Let us check the properties of the function $P_{2M}^{\ast}g$
where $g$ satisfies (\ref{gorthogonal1})-(\ref{gorthogonal3}). First of all,
it is clear from (\ref{gorthogonal1}) that $P_{2M}^{\ast}g\in X_{N}^{G}$ where
we have used the notation (\ref{XNG}).

By Green's formula (\ref{GreenGeneral}), applied for the operator $P_{2M}$ and
for $u=g$ we obtain
\[%
{\displaystyle\int_{G}}
P_{2M}^{\ast}g\cdot vdx=0\qquad\text{for all }v\text{ with }P_{2M}v=0
\]
which implies that the function $P_{2M}^{\ast}g$ satisfies $P_{2M}^{\ast
}g\perp X_{M}^{G}$ ($X_{M}^{G}$ defined as (\ref{XNG})).

By the general existence theorem for Elliptic BVP used already above (cf.
\cite{lions-magenes}, Theorem $5.3$, the Fredholmness property), we know that
a solution $g$ to problem (\ref{gorthogonal1})-(\ref{gorthogonal3}) exists for
those boundary data $\left\{  h_{j}\right\}  _{j=1}^{N-M}$ which satisfy only
a finite number of linear conditions (cf. \cite{lions-magenes}, conditions
(5.18)); these are determined by the solutions to the homogeneous adjoint
Elliptic BVP. Hence, it follows that the space $Y_{N-M}^{G}$ of the functions
$P_{2M}^{\ast}g$ where $g$ is a solution to (\ref{gorthogonal1}%
)-(\ref{gorthogonal3}) is infinite-dimensional.

\textbf{(5) }Let us construct a subspace of $X_{N}\setminus X_{M}$ which is
infinite-dimensional. We use the obvious inclusion $X_{N|G}\subset X_{N}^{G},$
$X_{M|G}\subset X_{M}^{G},$ where for a space of functions $Y\subset
L_{2}\left(  B\right)  $ the space $\ Y_{|G}$ consists of the restrictions of
the elements of $Y$ to the domain $G.$

First of all, we find an orthonormal basis $\left\{  v_{j}\right\}  _{j\geq1}$
in the infinite-dimensional space $Y_{N-M}^{G}$ (where the norm is $\left\Vert
\cdot\right\Vert _{L_{2}\left(  G\right)  }$ ); by the Gram-Schmidt
orthonormalization we obtain functions $g_{j}$ such that $v_{j}=P_{2M}^{\ast
}g_{j}$ for $j\geq1.$

Let us put $\varepsilon_{j}=\frac{1}{2^{j-1}}$ and use the density Theorem
\ref{Tdensity} to choose $u_{j}\in H^{2N}\left(  D\right)  $ with
\[
\left\Vert u_{j}-v_{j}\right\Vert _{L_{2}\left(  G\right)  }\leq
\varepsilon_{j}\qquad\text{for }j\geq1.
\]
The orthogonality of $v_{j}$ to $X_{M}^{G}$ infers $\operatorname*{dist}%
\left(  u_{j|G},X_{M}^{G}\right)  \geq1-\varepsilon_{j}$ in the $L_{2}\left(
G\right)  $ norm. Hence, $\operatorname*{dist}\left(  u_{j},X_{M}\right)
\geq1-\varepsilon_{j}$ in the $L_{2}\left(  D\right)  $ norm, hence
$u_{j}\notin X_{M}.$

Let us see that for every choice of the constants $\alpha_{j}$ holds
\[%
{\displaystyle\sum_{j=1}^{N-1}}
\alpha_{j}u_{j}\neq u_{N}.
\]
Indeed, by the triangle inequality for the norm $\left\Vert \cdot\right\Vert
_{L_{2}\left(  G\right)  }$ it follows
\begin{align*}
1+%
{\displaystyle\sum_{j=1}^{N-1}}
\left\vert \alpha_{j}\right\vert ^{2}  &  =\left\Vert v_{N}-%
{\displaystyle\sum_{j=1}^{N-1}}
\alpha_{j}v_{j}\right\Vert \\
&  =\left\Vert v_{N}-u_{N}+u_{N}-%
{\displaystyle\sum_{j=1}^{N-1}}
\alpha_{j}u_{j}+%
{\displaystyle\sum_{j=1}^{N-1}}
\alpha_{j}u_{j}-%
{\displaystyle\sum_{j=1}^{N-1}}
\alpha_{j}v_{j}\right\Vert \\
&  \leq\varepsilon_{N}+\left\Vert u_{N}-%
{\displaystyle\sum_{j=1}^{N-1}}
\alpha_{j}u_{j}\right\Vert +%
{\displaystyle\sum_{j=1}^{N-1}}
\left\vert \alpha_{j}\right\vert \varepsilon_{j}%
\end{align*}
or
\[
1-\varepsilon_{N}+%
{\displaystyle\sum_{j=1}^{N-1}}
\left(  \left\vert \alpha_{j}\right\vert ^{2}-\left\vert \alpha_{j}\right\vert
\varepsilon_{j}\right)  \leq\left\Vert u_{N}-%
{\displaystyle\sum_{j=1}^{N-1}}
\alpha_{j}u_{j}\right\Vert .
\]
Obviously
\[
1-\varepsilon_{N}+%
{\displaystyle\sum_{j=1}^{N-1}}
\left(  \frac{\varepsilon_{j}^{2}}{4}-\frac{\varepsilon_{j}}{2}\varepsilon
_{j}\right)  \leq1-\varepsilon_{N}+%
{\displaystyle\sum_{j=1}^{N-1}}
\left(  \left\vert \alpha_{j}\right\vert ^{2}-\left\vert \alpha_{j}\right\vert
\varepsilon_{j}\right)
\]
and since the left-hand side always exceed $1/4$, this ends the proof that the
system of functions $\left\{  u_{j|G}\right\}  _{j\geq1}$ is linearly
independent. Hence, the system $\left\{  u_{j}\right\}  _{j\geq1}$ is linearly
independent in the whole domain $D.$

As noted above $u_{j}\notin X_{M},$ hence $\operatorname*{span}\left\{
u_{j}\right\}  _{j\geq1}$ is the infinite-dimensional space we sought The
proof is finished.%

\endproof

We have the following prototype of Theorem \ref{TtransversalSpaces}, proved in
\cite{kounchevSozopol}.

\begin{corollary}
\label{Ctransversal}Let $M<N$ and $X_{M},$ $X_{N}$ satisfy Definition
\ref{Dhdimension} with
\[
\operatorname{hdim}\left(  X_{M}\right)  =M,\quad\operatorname{hdim}\left(
X_{N}\right)  =N.
\]
Assume that the differential operators $P_{2M}$ and $P_{2N}^{\prime},$
associated with $X_{M}$ and $X_{N},$ have trivial factorization operators by
the definition (\ref{P2Nproduct}), and trivial domain partitions $D=D_{1}$ and
$D=D_{1}^{\prime}$ by (\ref{domainPartition}). Then the space of solutions of
the Elliptic BVP (\ref{gorthogonal1})-(\ref{gorthogonal3}) where $G=D$ is a
subspace of the space
\[
Y=X_{N}\ominus X_{M}.
\]

\end{corollary}

The proof may be derived from the proof of Theorem \ref{TtransversalSpaces}
where we have put $G=D.$ Note that we \textbf{do not need} the $\left(
U\right)  _{s}$ condition for the operator $P_{2N}^{\prime}.$ Hence, strictly
speaking, Corollary \ref{Ctransversal} is not a special case of Theorem
\ref{TtransversalSpaces}.

Now we provide a generalization of Kolmogorov's notion of width from formula
(\ref{KolmogorovWidth}); without restricting the generality we assume that we
work only with symmetric subsets.

\begin{definition}
\label{Dwidth} Let $A$ be a centrally symmetric subset in $L_{2}\left(
B\right)  .$ For fixed integers $M\geq1$ and $N\geq0$ we define the
corresponding \textbf{Harmonic Width} by putting
\[
\operatorname*{hd}\nolimits_{M,N}\left(  K\right)  :=\inf_{X_{M},F_{N}%
}\operatorname*{dist}\left(  X_{M}%
{\textstyle\bigoplus}
F_{N},A\right)  ,
\]
where $\inf_{X_{M},F_{N}}$ is taken over all spaces $X_{M},F_{N}\subset
C^{\infty}\left(  B\right)  $ with
\begin{align*}
\operatorname{hdim}\left(  X_{M}\right)   &  =M\\
\dim\left(  F_{N}\right)   &  =N.
\end{align*}

\end{definition}

\section{Generalization of Kolmogorov's result about widths \label{Swidths}}

Next we prove results which are analogs to the original Kolmogorov's results
about widths in (\ref{KolmogorovWidth2}).

We denote by $F_{N}$ a \textbf{finite-dimensional} subspace of $L_{2}\left(
B\right)  $ of dimension $N.$ We denote the special subspaces for an elliptic
operator $P_{2p}=L_{2p}$ by
\begin{equation}
\widetilde{X}_{p}:=\left\{  u\in H^{2p}\left(  B\right)  :L_{2p}u\left(
x\right)  =0,\quad\text{for }x\in B\right\}  , \label{StildeM}%
\end{equation}
and the special finite-dimensional subspaces
\begin{equation}
\widetilde{F}_{N}:=\left\{  \psi_{j}:j\leq N\right\}  _{lin} \label{FtildeN}%
\end{equation}
where $\psi_{j}$ are the eigenfunctions from Theorem
\ref{TExpansionBerezanskii}.

\begin{theorem}
\label{TKolmogorovMultivariate} Let $K_{p}^{\ast}$ be the set defined in
(\ref{KpstarGENERAL}) as
\[
K_{p}^{\ast}:=\left\{  u\in H^{2p}\left(  B\right)  :%
{\displaystyle\int_{B}}
\left\vert L_{2p}u\left(  x\right)  \right\vert ^{2}dx\leq1\right\}  ,
\]
with a constant coefficient operator $L_{2p}$ which is uniformly strongly
elliptic in the domain $B.$ Let $X_{M}$ be a First Kind subspace of
$L_{2}\left(  B\right)  $ of Harmonic Dimension $M,$ according to Definition
\ref{Dhdimension}, i.e.
\[
\operatorname{hdim}\left(  X_{M}\right)  =M,
\]
and let $N\geq0$ be arbitrary.

1. If $M<p$ then
\[
\operatorname*{dist}\left(  X_{M}%
{\textstyle\bigoplus}
F_{N},K_{p}^{\ast}\right)  =\infty.
\]
Hence,
\[
\inf_{X_{M},F_{N}}\operatorname*{dist}\left(  X_{M}%
{\textstyle\bigoplus}
F_{N},K_{p}^{\ast}\right)  =\infty
\]
or equivalently,
\[
\operatorname*{hd}\nolimits_{M,N}\left(  K_{p}^{\ast}\right)  =\infty.
\]

2. If $M=p$ then
\[
\inf_{X_{p},F_{N}}\operatorname*{dist}\left(  X_{p}%
{\textstyle\bigoplus}
F_{N},K_{p}^{\ast}\right)  =\operatorname*{dist}\left(  \widetilde{X}_{p}%
{\textstyle\bigoplus}
\widetilde{F}_{N},K_{p}^{\ast}\right)  ,
\]
i.e.
\[
\operatorname*{hd}\nolimits_{p,N}\left(  K_{p}^{\ast}\right)
=\operatorname*{dist}\left(  \widetilde{X}_{p}%
{\textstyle\bigoplus}
\widetilde{F}_{N},K_{p}^{\ast}\right)  .
\]

\end{theorem}

\begin{remark}
In both cases we see that the special spaces $\widetilde{X}_{M}%
{\textstyle\bigoplus}
\widetilde{F}_{N}$ are extremizers among the large class of spaces $X_{M}%
{\textstyle\bigoplus}
F_{N}.$
\end{remark}

%

\proof
1. If we assume that $X_{M}$ and $\widetilde{X}_{p}$ are transversal the proof
is clear since $\widetilde{X}_{p}\subset K_{p}^{\ast}$ and there will be an
infinite-dimensional subspace in $\widetilde{X}_{p}\subset K_{p}^{\ast}$
containing at least one infinite axis with direction $f\in\widetilde{X}%
_{p}\setminus X_{M},$ such that
\[
\operatorname*{dist}\left(  X_{M}%
{\textstyle\bigoplus}
F_{N},f\right)  >0
\]
which implies
\[
\operatorname*{dist}\left(  X_{M}%
{\textstyle\bigoplus}
F_{N},K_{p}^{\ast}\right)  =\infty.
\]
If they are not transversal we remind that operators with analytic
coefficients satisfy the $\left(  U\right)  _{s}$ condition, and we may apply
Lemma \ref{LtransversalSpaces}.

2. For proving the second item, let us first note that $\widetilde{X}%
_{p}\subset X_{p}%
{\textstyle\bigoplus}
F_{N}.$ Indeed, since $\widetilde{X}_{p}\subset K_{p}^{\ast}$ the violation of
$\widetilde{X}_{p}\subset X_{p}%
{\textstyle\bigoplus}
F_{N}$ would imply that there exists an infinite axis $f$ in $K_{p}^{\ast}$
not contained in $X_{p}%
{\textstyle\bigoplus}
F_{N}$ which would immediately give
\[
\operatorname*{dist}\left(  X_{p}%
{\textstyle\bigoplus}
F_{N},K_{p}^{\ast}\right)  =\infty.
\]
Using the notations of Definition \ref{Dhdimension}, there exists a finite
cover $\bigcup D_{j}=B,$ and by Lemma \ref{LellipticCanonical} (applied for
$M=N=p$ ) it follows that on every subdomain $D_{j}$ holds $P_{2p}^{j}%
=C_{j}\left(  x\right)  L_{2p}$ for some function $C_{j}\left(  x\right)  .$
Thus we see that every $u\in X_{p}$ is a piecewise solution of $L_{2p}u=0$ on
$B,$ satisfying the interface conditions (\ref{Interfaces}) in Definition
\ref{Dhdimension}. Here we use an uniqueness theorem for "analytic
continuation" across the boundary argument (proved directly by Green's formula
(\ref{GreenGeneral}) as in \cite{okbook}, Lemma $20.10$ and the proof of
Theorem $20.11,$ p. $422$) that $u\in\widetilde{X}_{p},$ hence $X_{p}%
=\widetilde{X}_{p}.$

Further we follow the usual way as in \cite{lorentz} to see that
$\widetilde{F}_{N}$ is extremal among all finite-dimensional spaces $F_{N},$
i.e.
\[
\inf_{F_{N}}\operatorname*{dist}\left(  \widetilde{X}_{p}%
{\textstyle\bigoplus}
F_{N},K_{p}^{\ast}\right)  =\operatorname*{dist}\left(  \widetilde{X}_{p}%
{\textstyle\bigoplus}
\widetilde{F}_{N},K_{p}^{\ast}\right)  .
\]
This ends the proof.%

\endproof

We prove the following fundamental result which shows the mutual position of
two subspaces:

\begin{lemma}
\label{LtransversalSpaces}Assume the conditions of Theorem
\ref{TtransversalSpaces}. Let the integer $M_{1}\geq0.$ Then%
\[
\operatorname*{dist}\left(  X_{M}%
{\textstyle\bigoplus}
F_{M_{1}},X_{N}\right)  =\infty.
\]

\end{lemma}

The proof follows directly from Theorem \ref{TtransversalSpaces} since a
finite-dimensional subspace $F_{M_{1}}$ would not disturb the arguments there.

We obtain immediately the following result.

\begin{corollary}
Let us denote by $U_{N+1}$ the unit ball in $X_{N+1}$ in the $L_{2}\left(
B\right)  $ norm. Then
\[
\operatorname*{dist}\left(  X_{N},U_{N+1}\right)  =1.
\]

\end{corollary}

\begin{remark}
Lemma \ref{LtransversalSpaces} and especially the above Corollary may be
considered as a generalization in our setting of a theorem of Gohberg-Krein of
$1957$ (cf. \cite{lorentz}, Theorem $2$ on p. $137$ ) in a Hilbert space.
\end{remark}

We need the following intuitive result which is however not trivial.

\begin{lemma}
\label{LellipticCanonical}Let for the strongly elliptic differential operators
$L_{2N}=P_{2N}\left(  x;D_{x}\right)  $ and $P_{2M}=P_{2M}\left(
x;D_{x}\right)  $ of orders respectively $2N\leq2M$ in the domain $B,$ the
following inclusion hold
\[
X_{N}%
{\textstyle\bigcap}
H^{2M}\left(  B\right)  \subset X_{M}\setminus F,
\]
or
\begin{align*}
&  \left\{  u\in H^{2M}\left(  B\right)  :L_{2N}u\left(  x\right)  =0,\quad
x\in B\right\}  \subset\\
&  \subset\left\{  u\in H^{2M}\left(  B\right)  :P_{2M}u\left(  x\right)
=0,\quad x\in B\right\}  \setminus F,
\end{align*}
where $F\subset L_{2}\left(  B\right)  $ is a \textbf{finite-dimensional}
subspace of $L_{2}\left(  B\right)  .$ Then
\begin{equation}
P_{2M}\left(  x,D_{x}\right)  =P_{2M-2N}^{\prime}\left(  x,D_{x}\right)
L_{2N}\left(  x,D_{x}\right)  \label{P2N=cx}%
\end{equation}
for some strongly elliptic differential operator $P_{2M-2N}^{\prime}$ of order
$2M-2N.$
\end{lemma}

%

\proof
It is clear that the arguments for proving equality (\ref{P2N=cx}) are purely
local, and it suffices to consider only $x_{0}=0,$ or we assume that the
operator $L_{2N}$ has constant coefficients.

First, we assume that the polynomial $L_{2N}\left(  \zeta\right)  $ is
\emph{irreducible}. Then we consider the roots of the equation
\begin{equation}
L_{2N}\left(  \zeta\right)  =0\qquad\text{for }\zeta\in\mathbb{C}^{n}.
\label{L2Nzeta=0}%
\end{equation}
If $\zeta$ is a solution to (\ref{L2Nzeta=0}) then the function $v\left(
x\right)  =\exp\left(  \left\langle \zeta,x\right\rangle \right)  $ is a
solution to equation $L_{2N}v=0$ in the whole space. Hence
\[
P_{2M}v=P_{2M}\left(  x_{0};D_{x}\right)  v=P_{2M}\left(  x_{0};\zeta\right)
=0,
\]
and by a well-known result on division of polynomials in algebra \cite{walker}
(Theorem $9.7,$ p. $26$), the statement of the theorem follows.

Now let us assume that $L_{2N}$ is reducible and decomposed in two irreducible
factors $L_{2N}=Q_{2}Q_{1},$ which may be equal. Obviously, both polynomials
$Q_{1}$ and $Q_{2}$ are uniformly strongly elliptic. Since the solutions to
$Q_{1}u=0$ are also solutions to $L_{2N}$ it follows by the above that
\[
P_{2M}\left(  x,D_{x}\right)  =P_{2M-2N_{1}}^{\prime}\left(  x,D_{x}\right)
Q_{1}\left(  D_{x}\right)
\]
where $2N_{1}$ is the order of the operator $Q_{1}.$ Further, following the
standard arguments in \cite{lions-magenes}, by the uniform strong ellipticity
of the operator $Q_{1},$ for every $\zeta\in\mathbb{C}^{n},$ and for arbitrary
$s\geq2N_{1},$ there exists a solution $u\in H^{s}\left(  B\right)  $ to
equation
\[
Q_{1}u_{\zeta}\left(  x\right)  =e^{\left\langle \zeta,x\right\rangle }%
\qquad\text{for }x\in B.
\]
Let $\zeta\in\mathbb{C}^{n}$ be a solution to equation $Q_{2}\left(
\zeta\right)  =0.$ Obviously,
\[
L_{2N}u_{\zeta}=0
\]
hence, by the above it follows
\[
P_{2M}\left(  x,D_{x}\right)  u_{\zeta}=P_{2M-2N_{1}}^{\prime}\left(
x,D_{x}\right)  Q_{1}\left(  D_{x}\right)  u_{\zeta}=P_{2M-2N_{1}}^{\prime
}\left(  x,\zeta\right)  =0.
\]
It follows that $P_{2M-2N_{1}}^{\prime}\left(  x_{0},\zeta\right)  =0.$ We
proceed inductively if $L_{2N}$ has more than two irreducible factors.%

\endproof

\section{Second Kind spaces of Harmonic Dimension $N$ and widths
\label{SsecondKind}}

In order to make things more transparent, in Definition \ref{Dhdimension} we
avoided the maximal generality of the notions and considered only First Kind
spaces of Harmonic Dimension $N.$ Let us explain by analogy with the
one-dimensional case how do the "Second Kind" spaces of Harmonic Dimension $N$ appear.

In the one-dimensional case, if we have a finite-dimensional subspace
$X_{N}\subset C^{N}\left(  I\right)  $ then for a point $x_{0}\in I$ the
space
\[
Y:=\left\{  u\in X_{N}:u\left(  x_{0}\right)  =0\right\}
\]
is an $\left(  N-1\right)  -$dimensional subspace. We would like that our
notion of Harmonic Dimension $N$ behave in a similar way. For example, if
$X_{N}$ is defined as a set of solutions of an elliptic operator $P_{2N}$ by
\[
X_{N}:=\left\{  u\in H^{2N}\left(  B\right)  :P_{2N}u=0\quad\text{in
}B\right\}
\]
then it is natural to expect that the space
\[
Y:=\left\{  u\in X_{N}:u=0\quad\text{on }\partial B\right\}
\]
has Harmonic Dimension $N-1.$ A simple example is the space
\[
Y=\left\{  u\in H^{4}\left(  B\right)  :\Delta^{2}u=0\quad\text{in
}B,\ u=0\quad\text{on }\partial B\right\}  .
\]

On the other hand, it is Theorem \ref{TtransversalSpaces} and Corollary
\ref{Ctransversal} above which show that such Second Kind spaces of Harmonic
Dimension $N$ appear in a natural way when we consider the space
$X_{N}\circleddash X_{M}$ based on solutions of Elliptic BVP
(\ref{gorthogonal1})-(\ref{gorthogonal3}).

We give the following definition.

\begin{definition}
\label{DsecondKind} For an integer $M\geq1$ we say that the linear subspace
$X_{M}\subset L_{2}\left(  D\right)  $ is of \textbf{Second Kind} and has
\textbf{Harmonic Dimension} $M,$ and write
\[
\operatorname{hdim}\left(  X_{M}\right)  =M,
\]
if it satisfies all conditions of Definition \ref{Dhdimension} however with an
elliptic operator $P_{2N},$ with $N\geq M$ and all elements $u\in X_{M}$
satisfy $N-M$ boundary conditions
\[
B_{j}u=0\qquad\text{on }\partial D,\ j=1,2,...,N-M.
\]
Here the boundary operators $\left\{  B_{j}\right\}  _{j=1}^{N-M}$ are a
\textbf{normal system} of boundary operators defined on $\partial D,$ by
Definition \ref{Delliptic}, item 4).
\end{definition}

By a technique similar to the already used we may prove the following results
which generalize Theorem \ref{TKolmogorovMultivariate}. We assume that
$K_{p}^{\ast}$ is the set defined by (\ref{KpstarGENERAL}) with a strongly
elliptic \emph{constant coefficients} operator $L_{2p}.$ The space
$\widetilde{X}_{p}$ is defined by (\ref{StildeM}) and the space $\widetilde
{F}_{L}$ by (\ref{FtildeN}).

The following theorem is a generalization of item 1) in Theorem
\ref{TKolmogorovMultivariate}.

\begin{theorem}
\label{TSecondKindMlessp} Let $M<p$ and $L\geq0$ be arbitrary integer. Let
$X_{M}$ be a Second Kind space with Harmonic Dimension $N,$ i.e.
\[
\operatorname{hdim}\left(  X_{M}\right)  =M.
\]
Let $F_{L}$ be an $L-$dimensional subset of $L_{2}\left(  B\right)  .$ Then
\[
\operatorname*{dist}\left(  X_{M}%
{\textstyle\bigoplus}
F_{L},K_{p}^{\ast}\right)  =\infty.
\]

\end{theorem}

The proof of Theorem \ref{TSecondKindMlessp} follows with minor modifications
of Lemma \ref{LtransversalSpaces} (Theorem \ref{TtransversalSpaces}).

It is more non-trivial to consider the case $N=p.$ First we must prove the
following result.

\begin{lemma}
\label{LSecondKindM=p}Let $X_{p}$ be a Second Kind space of Harmonic Dimension
$p$ and $L\geq0$ be an arbitrary integer. Let $F_{L}$ be an $L-$dimensional
subset of $L_{2}\left(  B\right)  .$ Then
\[
\operatorname*{dist}\left(  X_{p}%
{\textstyle\bigoplus}
F_{L},K_{p}^{\ast}\right)  <\infty
\]
implies
\begin{equation}
\widetilde{X}_{p}\subset X_{p}. \label{XptildeConXp}%
\end{equation}
Let the elliptic operator $P_{2M}$ and the boundary operators $\left\{
B_{j}\right\}  _{j=1}^{M-p}$ be associated with $X_{p}$ by Definition
\ref{DsecondKind}. Then (\ref{XptildeConXp}) implies the following
factorizations:
\begin{align*}
P_{2M}  &  =P_{2M-2p}^{\prime}L_{2p}\\
B_{j}  &  =B_{j}^{\prime}L_{2p}\qquad\text{for }j=1,2,...,M-p.
\end{align*}
The operator $P_{2M-2p}^{\prime}$ is uniformly strongly elliptic in $D,$ and
the boundary operators $\left\{  B_{j}^{\prime}\right\}  _{j=1}^{M-p}$ form a
normal system which covers the operator $P_{2M-2p}^{\prime}.$
\end{lemma}

Finally, the following generalization of item 2) in Theorem
\ref{TKolmogorovMultivariate} may be proved. It shows that one needs to take
into account the index of the Elliptic BVP involved.

\begin{theorem}
\label{TKolmogorovSecondKind} Let us consider those spaces $X_{p}$ of Second
Kind with Harmonic Dimension $p$ for which
\[
\operatorname*{dist}\left(  X_{p}%
{\textstyle\bigoplus}
F_{L},K_{p}^{\ast}\right)  <\infty
\]
with associated operators $P_{2M}$ and boundary operators $\left\{
B_{j}\right\}  _{j=1}^{M-p}.$ Following the notations of Lemma
\ref{LSecondKindM=p}, let us denote by $\mathcal{N}$ the following space of
solutions $w\in H^{2M-2p}\left(  D\right)  $ of the Elliptic BVP on the domain
$D$:%
\begin{align*}
P_{2M-2p}^{\prime}w  &  =0\qquad\text{on }D\\
B_{j}^{\prime}w  &  =0,\ \text{on }\partial D,\ \text{for }j=1,2,...,M-p
\end{align*}
Then the following equality holds \
\begin{equation}
\inf_{X_{p},F_{L}}\left\{  \operatorname*{dist}\left(  X_{p}%
{\textstyle\bigoplus}
F_{L},K_{p}^{\ast}\right)  :\text{ }\dim\left(  \mathcal{N}\right)
+L=L_{1}\right\}  =\operatorname*{dist}\left(  \widetilde{X}_{p}%
{\textstyle\bigoplus}
\widetilde{F}_{L_{1}},K_{p}^{\ast}\right)  .\nonumber
\end{equation}

\end{theorem}

From the theory of Elliptic BVP is known that $\dim\left(  \mathcal{N}\right)
<\infty$ (cf. \cite{lions-magenes}, Theorem $5.3,$ chapter $2,$ section
$5.3$). Let us denote by $\left\{  w_{s}\right\}  _{s=1}^{\dim\left(
\mathcal{N}\right)  }$ a basis of the space $\mathcal{N},$ and by $u_{s}$ a
fixed solution to $L_{2p}u_{s}=w_{s}.$ The main point in the proof of Theorem
\ref{TKolmogorovSecondKind} is that arbitrary solution $u$ to equation
$P_{2M}u=0$ may be expressed as
\[
u=%
{\displaystyle\sum_{s=1}^{\dim\left(  \mathcal{N}\right)  }}
\lambda_{s}u_{s}+v
\]
where $v$ is a solution to $L_{2p}v=0.$

\section{Appendix, Proof of Theorem \ref{TExpansionBerezanskii}
\label{Sappendix}}

%

\proof
\textbf{(1)} We consider the following auxiliary elliptic \emph{eigenvalue
problem}
\begin{align}
L_{2p}L_{2p}^{\ast}\phi\left(  x\right)   &  =\lambda\phi\left(  x\right)
\qquad\qquad\text{on }B,\label{ffiBVP1}\\
B_{j}\phi\left(  y\right)   &  =S_{j}\phi\left(  y\right)  =0\qquad\text{for
}j=0,1,...,p-1,\text{ for }y\in\partial B. \label{ffiBVP2}%
\end{align}
Since this is the Dirichlet problem for the operator $L_{2p}^{\ast}L_{2p}$ it
is a classical fact that (\ref{ffiBVP1})-(\ref{ffiBVP2}) is a \emph{regular
Elliptic BVP} considered in the Sobolev space $H^{2p}\left(  B\right)  ,$ as
defined in Definition \ref{Delliptic}. Also, it is a classical fact that the
Dirichlet problem is a self-adjoint problem (cf. \cite{lions-magenes}, Remark
$2.4$ in section $2.4$ and Remark $2.6$ in section $2.5$, chapter $2$).

Hence, we may apply the main results about the Spectral theory of regular
self-adjoint Elliptic BVP. We refer to \cite{egorov-shubin0} (section $3$ in
chapter $2,$ p. $122,$ Theorem $2.52$) and to references therein.

By the uniqueness Lemma \ref{LuniquenessDirichlet} the eigenvalue problem
(\ref{ffiBVP1})-(\ref{ffiBVP2}) has only zero solution for $\lambda=0.$ It has
eigenfunctions $\phi_{k}\in H^{2p}\left(  B\right)  $ with eigenvalues
$\lambda_{k}>0$ for $k=1,2,3,...$ for which $\lambda_{k}\longrightarrow\infty$
as $k\longrightarrow\infty.$

\textbf{(2)} Next, in the Sobolev space $H^{2p}\left(  B\right)  ,$ we
consider the problem:
\begin{align}
L_{2p}L_{2p}^{\ast}\varphi\left(  x\right)   &  =\phi_{k}\left(  x\right)
\qquad\qquad\text{on }B\label{fiBVP1}\\
B_{j}\varphi\left(  y\right)   &  =S_{j}\varphi\left(  y\right)
=0\qquad\text{for }j=0,1,...,p-1,\text{ for }y\in\partial B. \label{fiBVP2}%
\end{align}
Obviously, the Elliptic BVP defined by problem (\ref{fiBVP1})-(\ref{fiBVP2})
coincides with the Elliptic BVP defined by (\ref{ffiBVP1})-(\ref{ffiBVP2}) up
to the right-hand sides, and all remarks there hold as well. Hence, problem
(\ref{fiBVP1})-(\ref{fiBVP2}) has\emph{ unique solution} $\varphi_{k}\in
H^{2p}\left(  B\right)  .$ We put
\[
\psi_{k}=L_{2p}^{\ast}\varphi_{k}.
\]
Hence, $L_{2p}\psi_{k}=\phi_{k}.$ We infer that on the boundary $\partial B$
hold the equalities $B_{j}L_{2p}\psi_{k}=B_{j}\phi_{k}$ and $S_{j}L_{2p}%
\psi_{k}=S_{j}\phi_{k};$ since $\phi_{k}$ are solutions to (\ref{ffiBVP1}%
)-(\ref{ffiBVP2}) it follows
\begin{equation}
B_{j}L_{2p}\psi_{k}\left(  y\right)  =S_{j}L_{2p}\psi_{k}\left(  y\right)
=0\qquad\text{for }j=0,1,...,p-1,\text{ for }y\in\partial B. \label{psiBVP}%
\end{equation}

We will prove that $\psi_{k}$ are solutions to problem (\ref{eigen1Multi}%
)-(\ref{eigen2Multi}), they are mutually \textbf{orthogonal,} and they are
also orthogonal to the space $\left\{  v\in H^{2p}:L_{2p}v=0\right\}  .$

\textbf{(3) }Let us see that
\[
L_{2p}^{\ast}L_{2p}\psi_{k}=\lambda_{k}\psi_{k}.
\]
By the definition of $\psi_{k}$ this is equivalent to
\[
L_{2p}^{\ast}L_{2p}L_{2p}^{\ast}\varphi_{k}=\lambda_{k}L_{2p}^{\ast}%
\varphi_{k};
\]
from $L_{2p}L_{2p}^{\ast}\varphi_{k}=\phi_{k}$ this is equivalent to
\[
L_{2p}^{\ast}\phi_{k}=\lambda_{k}L_{2p}^{\ast}\varphi_{k}%
\]
On the other hand, by the basic properties of $\phi_{k}$ and $\varphi_{k},$ we
have obviously $L_{2p}L_{2p}^{\ast}\phi_{k}=\lambda_{k}L_{2p}L_{2p}^{\ast
}\varphi_{k},$ hence
\[
L_{2p}L_{2p}^{\ast}\left(  \phi_{k}-\lambda_{k}\varphi_{k}\right)  =0.
\]
Note that both $\phi_{k}$ and $\varphi_{k}$ satisfy the same zero Dirichlet
boundary conditions, namely (\ref{ffiBVP2}) and (\ref{fiBVP2}). Hence, by the
uniqueness Lemma \ref{LuniquenessDirichlet} it follows that $\phi_{k}%
-\lambda_{k}\varphi_{k}=0$ which implies $L_{2p}^{\ast}L_{2p}\psi_{k}%
=\lambda_{k}\psi_{k}.$ Thus we see that $\psi_{k}$ is a solution to problem
(\ref{eigen1Multi})-(\ref{eigen2Multi}) and does not satisfy $L_{2p}\psi=0.$

\textbf{(4)} The orthogonality to the subspace $\left\{  v\in H^{2p}%
:L_{2p}v=0\right\}  $ follows easily from the Green formula
(\ref{GreenGeneral}) applied to the operator $L_{2p}^{\ast}L_{2p},$
\begin{align*}
&
{\displaystyle\int_{D}}
\left(  L_{2p}^{\ast}L_{2p}\psi_{k}\cdot v-L_{2p}\psi_{k}\cdot L_{2p}v\right)
dx\\
&  =%
{\displaystyle\sum_{j=0}^{2p-1}}
{\displaystyle\int_{\partial D}}
\left(  S_{j}L_{2p}\psi_{k}\cdot C_{j}v-B_{j}L_{2p}\psi_{k}\cdot
T_{j}v\right)
\end{align*}
in which substitute the zero boundary conditions (\ref{psiBVP}) of $\psi_{k},$
and equality
\[%
{\displaystyle\int_{D}}
L_{2p}^{\ast}L_{2p}\psi_{k}\cdot vdx=\lambda_{k}%
{\displaystyle\int_{D}}
\psi_{k}\cdot vdx.
\]

The orthonormality of the system $\left\{  \psi_{k}\right\}  _{k=1}^{\infty}$
follows now easily by the equality
\[
\lambda_{k}%
{\displaystyle\int}
\psi_{k}\psi_{j}dx=%
{\displaystyle\int}
L_{2p}^{\ast}L_{2p}\psi_{k}\psi_{j}dx=%
{\displaystyle\int}
L_{2p}\psi_{k}L_{2p}\psi_{j}dx=%
{\displaystyle\int}
\phi_{k}\phi_{j}dx
\]
and the orthogonality of the system $\left\{  \phi_{k}\right\}  _{k=1}%
^{\infty}.$

\textbf{(5)} For the completeness of the system $\left\{  \psi_{k}\right\}
_{k=1}^{\infty}$, let us assume that for some $f\in L_{2}\left(  B\right)  $
holds
\begin{equation}%
{\displaystyle\int_{B}}
f\cdot\psi_{k}dx=%
{\displaystyle\int_{B}}
f\cdot\psi_{k}^{\prime}dx=0\qquad\text{for all }k\geq1. \label{fOrthogonal}%
\end{equation}
Then the Green formula (\ref{GreenGeneral}) implies
\begin{align*}
0  &  =\lambda_{k}%
{\displaystyle\int_{B}}
f\cdot\psi_{k}dx=%
{\displaystyle\int_{B}}
f\cdot L_{2p}^{\ast}L_{2p}\psi_{k}dx=%
{\displaystyle\int_{B}}
L_{2p}f\cdot L_{2p}\psi_{k}dx\\
&  =%
{\displaystyle\int_{B}}
L_{2p}f\cdot\phi_{k}dx\qquad\text{for all }k\geq1.
\end{align*}
By the completeness of the system $\left\{  \phi_{k}\right\}  _{k\geq1}$ this
implies that $L_{2p}f=0.$ From the second orthogonality in (\ref{fOrthogonal})
follows that $f\equiv0,$ and this ends the proof of the completeness of the
system $\left\{  \psi_{j}^{\prime}\right\}  _{j=1}^{\infty}\bigcup\left\{
\psi_{j}\right\}  _{j=1}^{\infty}.$%

\endproof

We have used above the following simple result.

\begin{lemma}
\label{LuniquenessDirichlet} The solution to problem (\ref{ffiBVP1}%
)-(\ref{ffiBVP2}) for $\lambda=0$ is unique.
\end{lemma}

%

\proof
From Green's formula (\ref{GreenGeneral}) we obtain
\[%
{\displaystyle\int_{B}}
\left[  L_{2p}\phi\right]  ^{2}dx-%
{\displaystyle\int}
\phi\cdot L_{2p}^{\ast}L_{2p}\phi dx=%
{\displaystyle\sum_{j=1}^{p}}
{\displaystyle\int_{\partial B}}
\left(  S_{j}\phi\cdot C_{j}L_{2p}\phi-B_{j}\phi\cdot T_{j}L_{2p}\phi\right)
d\sigma_{y},
\]
hence $L_{2p}\phi=0.$

Now for arbitrary $v\in H^{2p}\left(  B\right)  $ by the same Green's formula
we obtain
\[%
{\displaystyle\int_{B}}
\left(  L_{2p}\phi\cdot v-\phi\cdot L_{2p}^{\ast}v\right)  dx=%
{\displaystyle\sum_{j=1}^{p}}
{\displaystyle\int_{\partial B}}
\left(  S_{j}\phi\cdot C_{j}v-B_{j}\phi\cdot T_{j}v\right)  d\sigma_{y}=0,
\]
hence
\[%
{\displaystyle\int_{B}}
\phi\cdot L_{2p}^{\ast}vdx=0.
\]
From the local existence theorem for elliptic operators (cf.
\cite{lions-magenes}) it follows that for arbitrary $f\in L_{2}\left(
B\right)  $ we may solve the elliptic equation $L_{2p}^{\ast}v=f$ with $v\in
H^{2p}\left(  B\right)  .$ From the density of $H^{2p}\left(  B\right)  $ in
$L_{2}\left(  B\right)  $ we infer $\phi\equiv0.$

This ends the proof.%

\endproof

\section{Conclusion and open problems}

As in Approximation, Spline and Wavelet Theory (\cite{kounchev1992TAMS},
\cite{okbook}), in the present research solutions of higher order elliptic
equation have shown flexibility which enabled a natural multidimensional
generalization of Kolmogorov's theory of widths with successful application to
multidimensional sets $K_{p}^{\ast}.$ Also, new features of Jackson type
theorems have been disclosed in Theorem \ref{TKolmogorovMultivariate}, which
shows that one needs components of different dimensions: $X_{p}$ and $F_{L}$
are of different types.

It may come as a big surprise, but the present research shows this
unambiguously, that in many issues one has to give up the convenient
simplistic understanding of the multidimensional case, in particular by
realizing that the finite-dimensional subspaces in $C^{N}\left(  D\right)  ,$
for domains $D\subset\mathbb{R}^{n}$ for $n\geq2,$ \textbf{do not} serve the
same job as the finite-dimensional subspaces in $C^{N}\left(  D\right)  $ for
intervals $D\subset\mathbb{R}^{1},$ and one has to replace them by a lot more
sophisticated objects, namely by the spaces having Harmonic Dimension $N.$
This is the \textbf{main conclusion} of the present research based on the
successful application of the new Harmonic Widths to explaining the structure
of the sets $K_{p}^{\ast}.$

Beyond the motivational problems mentioned in the Introduction, one may
formulate several other \textbf{open problems}:

\begin{enumerate}
\item First of all, one has to study basic questions about the sets having
Harmonic Dimension, by considering the sets $X_{M}\bigcap X_{N},$ $X_{M}%
{\textstyle\bigoplus}
X_{N},$ $X_{M}%
{\textstyle\bigotimes}
X_{N},$ and finding their Harmonic Dimension (if it exists!), etc.

\item Secondly, one has to prove a generalization of a theorem of S. Bernstein
about differentiable Markov systems (or, in the case of differentiability,
Extended Complete Chebyshev systems, in the terminology of
\cite{karlinstudden}). As remarked in \cite{kreinnudelman} (after the proof of
Theorem $4.2$ in chapter $2$), S. Bernstein dealt with even stronger
statement, namely, he was seeking Descartes systems (cf. \cite{kreinnudelman}%
). This needs the factorization of elliptic PDOs into $N$ elliptic operators
of second order. These operators will be obviously pseudo-differential,
\cite{hoermanderVol3}.

\item In this context, one has to check that the maximal generality of the
theory in the present paper will be achieved by considering elliptic
pseudo-differential operators.

\item New Jackson type theorems are suggested by the widths reasons:\ the
simplest way to state them is to consider spaces defined by
\[
\left\{  u:\left\vert L_{2p}u\left(  x\right)  \right\vert \leq1\qquad
\text{for }x\in D\right\}  .
\]
By arguments similar to the proof of Theorem \ref{TKolmogorovMultivariate} one
is convinced that a reasonable Jackson type theorem may be proved only for
operators $P_{2N}$ of the form $P_{2N}=P_{2N-2p}^{\prime}L_{2p},$ i.e. one has
to approximate through functions $u_{N}$ in the spaces $\left\{
u:P_{2N}u=0\quad\text{in }D\right\}  .$ In the case of polyharmonic operator
Jackson type results have been proved in \cite{kounchev1991Hanstholm}.

\item One has to find a proper discrete version of the present research which
will be essential for the applications to Compressed Sensing, compare the role
of Gelfand's widths in \cite{donoho}, \cite{devore}.

\item Although we have mentioned the Chebyshev systems in passing, an
important point of the present research is the generalization of the Extended
Complete Chebyshev systems of order $N$ (discussed in more detail in
\cite{kounchev2008-ChebyshevSystems}) which is the ground for the spaces
having Harmonic Dimension $N.$ One has to specify more precisely which are the
elliptic differential/pseudodifferential operators acceptable for a
Multidimensional Chebyshev system. This has to be considered in the context of
S. Bernstein's one-dimensional result, mentioned in Proposition
\ref{PstructureMarkovSystem}.

\item In the same direction, let us recall that one-dimensional Chebyshev
systems are important for the qualitative theory of ODEs, in particular for
\emph{Sturmian type} of theorems, cf. e.g. \cite{arnold}, \cite{arnoldProblem}%
. There has been a long search for proper multidimensional generalizations of
Chebyshev systems. The standard generalization by means of zero set property
fails to produce a non-trivial multidimensional system and this is the content
of the theorem of Mairhuber, cf. the thorough discussion in
\cite{kreinnudelman} (chapter $2,$ section $1.1$). In general, zero set
properties and intersections are not a reliable reference point for
multidimensional Analysis. Indeed, let us recall that polyharmonic (and even
harmonic) functions do not have simple zero sets, however they are solutions
to nice Dirichlet problems (\ref{DirichletMultivariate}) and for that reason
are considered to be a genuine Multidimensional Chebyshev system as we have
defined it in (\ref{S2M}).

V.I. Arnold discusses the importance of the Chebyshev systems in his Toronto
lectures, June 1997, Lecture 3: Topological Problems in Wave Propagation
Theory and Topological Economy Principle in Algebraic Geometry. Fields
Institute Communications, available online at\emph{
http://www.pdmi.ras.ru/\symbol{126}arnsem/Arnold/arn-papers.html. } On p. 8 he
writes that \textquotedblright Even the Sturm theory is missing in higher
dimensions. This is an interesting phenomenon. All attempts that I know to
extend Sturm theory to higher dimensions failed. For instance, you can find
such an attempt in the Courant-Hilbert's book, in chapter 6, but it is wrong.
The topological theorems about zeros of linear combinations for higher
dimensions, which are attributed there to Herman, are wrong even for the
standard spherical Laplacian.\textquotedblright\ The attempts to mimic
multivariate Chebyshev systems are present in the works of V.I. Arnold in the
context of multivariate Sturm type of theorems, see in particular problem
1996-5 in \cite{arnoldProblem}. In view of these efforts of V.I. Arnold, one
might try to apply the present framework for obtaining multidimensional Sturm
type theorems.
\end{enumerate}

\end{document}